\documentclass{IEEEtran}
\usepackage[T1]{fontenc}
\usepackage[latin9]{inputenc}
\usepackage{color}
\usepackage{array}
\usepackage{float}
\usepackage{mathrsfs}
\usepackage{multirow}
\usepackage{amsmath}
\usepackage{amssymb}
\usepackage{graphicx}
\usepackage{setspace}
\usepackage{esint}

\makeatletter

\providecommand{\tabularnewline}{\\}
\floatstyle{ruled}
\newfloat{algorithm}{tbp}{loa}
\providecommand{\algorithmname}{Algorithm}
\floatname{algorithm}{\protect\algorithmname}

\usepackage{url}
\usepackage{subfigure}
\usepackage{epsfig}
\usepackage{subfigmat}
\usepackage{amsthm}\usepackage{tensor}\usepackage{amsfonts}
\usepackage{mathrsfs}

\newcounter{assump}
\newcounter{defn}

\theoremstyle{plain}
\newtheorem{theorem}{\bf Theorem}
\newtheorem{lemma}[theorem]{\bf Lemma}
\newtheorem{proposition}[theorem]{\bf Proposition}

\theoremstyle{remark}
\newtheorem{assumption}[assump]{\bf Assumption}
\newtheorem{definition}[defn]{\bf Definition}
\newtheorem{remark}[defn]{\bf Remark}

\makeatother

\begin{document}

\title{Distributed Estimation using\\Bayesian Consensus Filtering}

\author{Saptarshi~Bandyopadhyay~\IEEEmembership{Student~Member,~IEEE},
Soon-Jo~Chung~\IEEEmembership{Senior~Member,~IEEE} %
\thanks{The authors are with the Department of Aerospace Engineering and Coordinated
Science Laboratory, University of Illinois at Urbana-Champaign, Urbana,
IL 61801, USA. (email: bandyop2@illinois.edu; sjchung@illinois.edu) 

This research was supported by AFOSR grant FA95501210193.%
} }
\maketitle
\begin{abstract}
We present the Bayesian consensus filter (BCF) for tracking a moving
target using a networked group of sensing agents and achieving consensus
on the best estimate of the probability distributions of the target's
states. Our BCF framework can incorporate nonlinear target dynamic
models, heterogeneous nonlinear measurement models, non-Gaussian uncertainties,
and higher-order moments of the locally estimated posterior probability
distribution of the target's states obtained using Bayesian filters.
If the agents combine their estimated posterior probability distributions
using a logarithmic opinion pool, then the sum of Kullback--Leibler
divergences between the consensual probability distribution and the
local posterior probability distributions is minimized. Rigorous stability
and convergence results for the proposed BCF algorithm with single
or multiple consensus loops are presented. Communication of probability
distributions and computational methods for implementing the BCF algorithm
are discussed along with a numerical example.
\end{abstract}

\section{Introduction}

In this paper, the term \textit{consensus} means reaching an agreement
across the network, regarding a certain subject of interest called
the \textit{target} dynamics. Distributed and networked groups of
agents can sense the target, broadcast the acquired information, and
reach an agreement on the gathered information using consensus algorithms.
Potential applications of distributed estimation tasks include environment
and pollution monitoring, tracking dust or volcanic ash clouds, tracking
mobile targets such as flying objects or space debris using distributed
sensor networks, etc. Consensus algorithms are extensively studied
in controls \nocite{Ref:Saber04,Ref:Jadbabaie03,Ref:Shah06,Ref:Boyd04,Ref:Ren05TAC,Ref:Saber07}
\cite{Ref:Saber04}--\cite{Ref:Saber07}, distributed optimization
\nocite{Ref:Tsitsiklis05,Ref:Tsitsiklis86,Ref:Martinez10,Ref:Nedic09}
\cite{Ref:Tsitsiklis05}--\cite{Ref:Nedic09}, and distributed estimation
problems \nocite{Ref:Borkar82,Ref:Moura11,Ref:Giannakis09,Ref:Murray05,Ref:Coates04}
\cite{Ref:Borkar82}--\cite{Ref:Coates04}. Strictly speaking, \textit{consensus}
is different from the term \textit{distributed estimation}, which
refers to finding the best estimate of the target, given a distributed
network of sensing agents. 

Many existing algorithms for distributed estimation \nocite{Ref:Borkar82,Ref:Moura11,Ref:Leonard10,Ref:Ahmed13,Ref:Demetriou11,Ref:Moura10,Ref:Hadaegh07}
\cite{Ref:Borkar82}--\cite{Ref:Hadaegh07} aim to obtain the estimated
mean (first moment of the estimated probability distribution) of the
target dynamics across the network, but cannot incorporate nonlinear
target dynamics, heterogeneous nonlinear measurement models, non-Gaussian
uncertainties, or higher-order moments of the locally estimated posterior
probability distribution of the target's states. It is difficult to
recursively combine local mean and covariance estimates using a linear
consensus algorithm because the dimension of the vector transmitted
by each agent increases linearly with time due to correlated process
noise \cite{Ref:Gupta13} and the covariance update equation is usually
approximated by a consensus gain \cite{Ref:Saber09}.

\begin{figure}
\begin{centering}
\includegraphics[bb=510bp 100bp 1000bp 560bp,clip,width=2.8in]{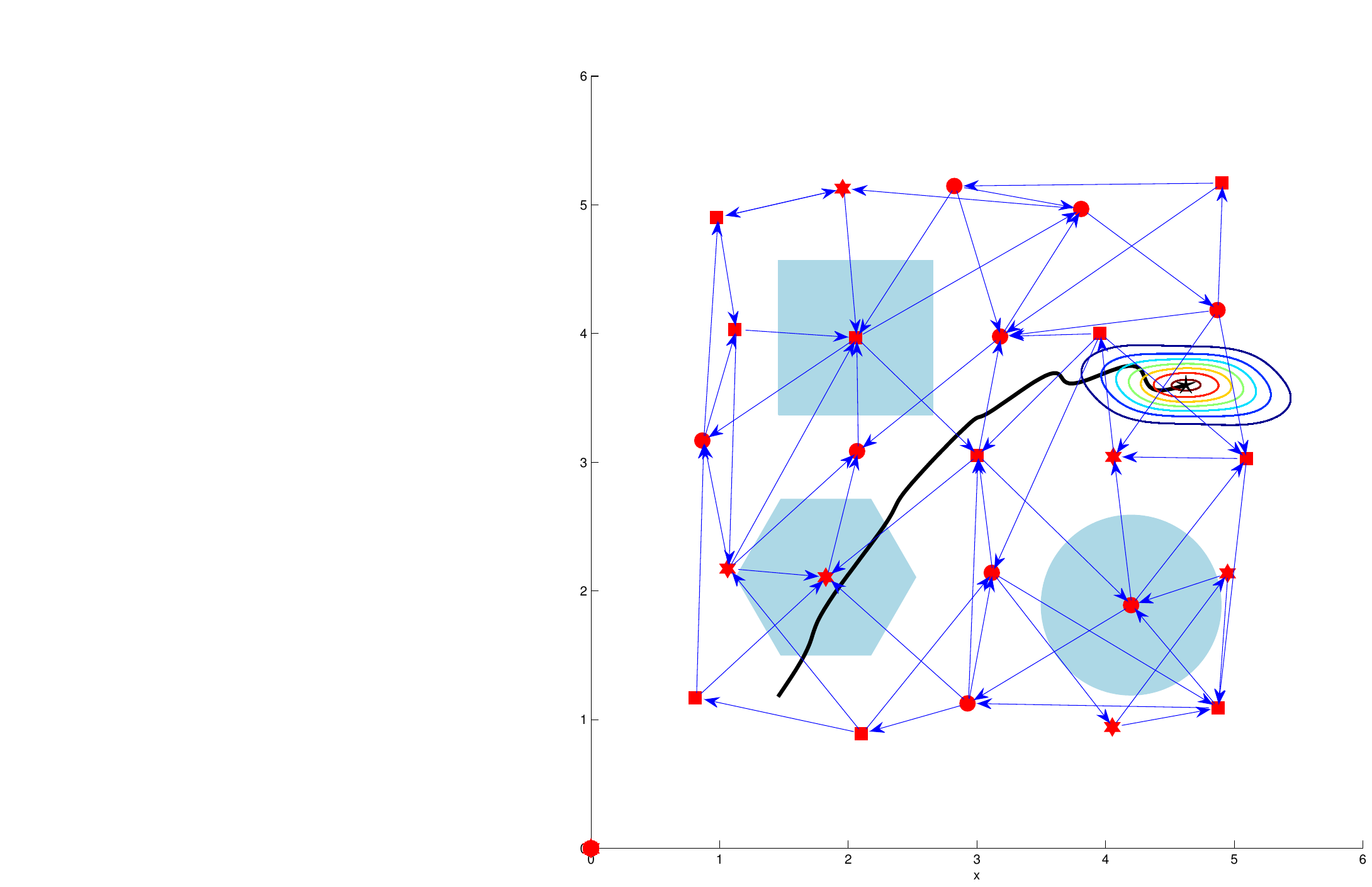}
\par\end{centering}

\caption{A heterogeneous sensor network tracks a target, whose current position
is marked by {\large{$\star$}}. The sensing region of some sensors
are shown in light blue. The sensors use the Hierarchical BCF algorithm
to estimate the probability distribution of the target's position
and reach a consensus across the network. \label{fig:BCF-overview}}

\vspace{-10pt}
\end{figure}

Multi-agent tracking or sensing networks are deployed in a distributed
fashion when the target dynamics have complex temporal and spatial
variations. Hence, it is necessary to preserve the complete information
captured in the locally estimated posterior probability distribution
of the target's states while achieving consensus across the network.
For example, while tracking an orbital debris in space, the uncertainty
in position along the direction of velocity is much larger than that
orthogonal to the velocity direction, and this extra information is
lost if a linear consensus algorithm is used to combine the estimated
mean from multiple tracking stations. As shown in Fig \ref{fig:BCF-overview},
the contour plot represents the consensual probability distribution
of the target's final position, where the uncertainty along the velocity
direction is relatively larger than that in the orthogonal direction.

The main objective of this paper is to extend the scope of distributed
estimation algorithms to track targets with general nonlinear dynamic
models with stochastic uncertainties, thereby addressing the aforementioned
shortcomings. Bayesian filters \nocite{Ref:DeGroot75,Ref:Jeffreys61,Ref:Subrahmaniam79,Ref:Pearl88}
\cite{Ref:DeGroot75}--\cite{Ref:Pearl88} recursively calculate the
probability density/mass function of the beliefs and update them based
on new measurements. The main advantage of Bayesian filters over Kalman
filter--based methods \cite{Ref:Kalman60,Ref:Julier04} for estimation
of nonlinear target dynamic models is that no approximation is needed
during the filtering process. In other words, the complete information
about the dynamics and uncertainties of the model can be incorporated
in the filtering algorithm. However, Bayesian filtering is computationally
expensive. Advances in computational capability have facilitated the
implementation of Bayesian filters for robotic localization and mapping
\nocite{Ref:Thrun00,Ref:Burgard96,Ref:Diard03,Ref:Thrun05} \cite{Ref:Thrun00}--\cite{Ref:Thrun05}
as well as planning and control \nocite{Ref:Boutilier99,Ref:Kaelbling98,Ref:Punska99}
\cite{Ref:Boutilier99}--\cite{Ref:Punska99}. Practical implementation
of these algorithms, in their most general form, is achieved using
particle filtering \cite{Ref:Pearl88,Ref:Arulampalam02} and Bayesian
programming \cite{Ref:Lebeltel04,Ref:Chen05}. This paper focuses
on developing a consensus framework for distributed Bayesian filters. 

The \textit{statistics} literature deals with the problem of reaching
a consensus among individuals in a complete graph, where each individual's
opinion is represented as a probability distribution \cite{Ref:DeGroot74,Ref:Genest86};
and under select conditions, it is shown that consensus is achieved
within the group \nocite{Ref:Bacharach79,Ref:Chatterjee77,Ref:French81}
\cite{Ref:Bacharach79}--\cite{Ref:French81}. Exchange of beliefs
in decentralized systems, under communication constraints, is considered
in \cite{Ref:Velagapudi07,Ref:Yuksel09}. Algorithms for combining
probability distributions within the exponential family, i.e., a limited
class of unimodal distributions that can be expressed as an exponential
function, are presented in \cite{Ref:Fraser12,Ref:Hlinka12}. If the
target's states are discrete random variables, then the local estimates
can be combined using a tree-search algorithm \cite{Ref:Pavlin10}
or a linear consensus algorithm \cite{Ref:Jadbabaie12}. In contrast,
this paper focuses on developing generalized Bayesian consensus algorithms
with rigorous convergence analysis for achieving consensus across
the network without any assumption on the shape of local prior or
posterior probability distributions. The proposed distributed estimation
using Bayesian consensus filtering aims to reach an agreement across
the network on the best estimate, in the information theoretic sense,
of the probability distribution of the target's states.

\subsection{Paper Contributions and Organization \label{sub:Paper-Contributions-and}}

\begin{figure}
\begin{centering}
\includegraphics[bb=100bp 270bp 770bp 1060bp,clip,width=3.2in]{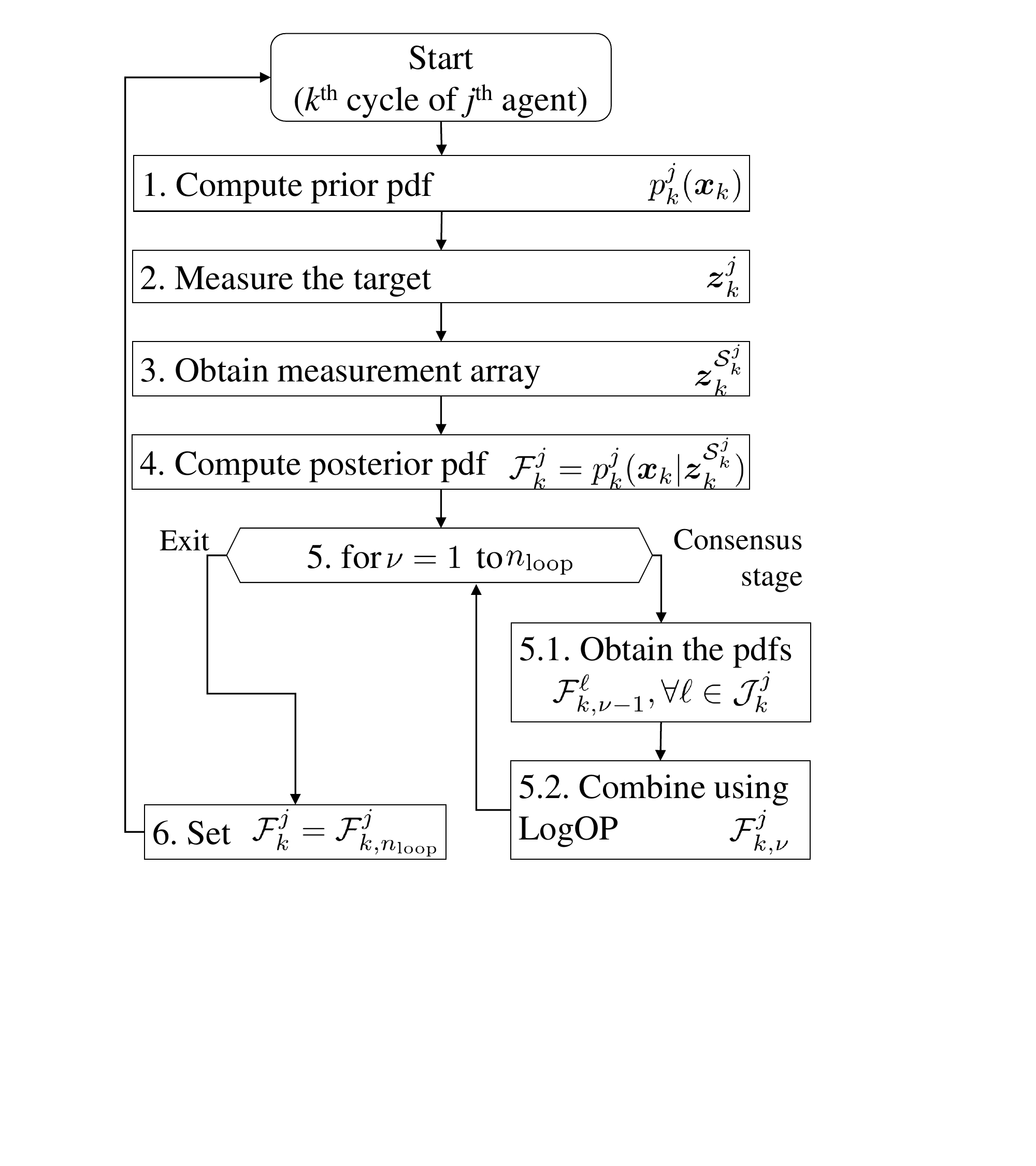}
\par\end{centering}

\caption{Flowchart for BCF--LogOP algorithm describing the key steps for a
single agent in a single time step. Steps 1--4 represent the Bayesian
filtering stage, while step 5 represents the consensus stage. \label{fig:Flowchart-BCF-LogOP}}

\vspace{-10pt}
\end{figure}

In this paper, we assume that agents generate their local estimate
of the posterior probability distribution of the target's states using
Bayesian filters with/without measurement exchange with neighbors.
Then, we develop algorithms for combining these local estimates, using
the logarithmic opinion pool (LogOP), to generate the consensual estimate
of the probability distribution of the target's states across the
network. Finally, we introduce the Bayesian consensus filter (BCF),
where the local prior estimates of the target\textquoteright{}s states
are first updated and the local posterior probability distributions
are recursively combined during the consensus stage, so that the agents
can estimate the consensual probability distribution of the target's
states while simultaneously maintaining consensus across the network.
The flowchart for the algorithm is shown in Fig. \ref{fig:Flowchart-BCF-LogOP}
and its pseudo-code is given in \textbf{Algorithm} \textbf{\ref{alg:Dynamic-Bayesian-Consensus}}. 

The first contribution of this paper is the LogOP--based consensus
algorithm for combining posterior probability distributions during
the consensus stage and achieving consensus across the network. As
discussed in Section \ref{sub:Logarithmic-Opinion-Pool}, this is
achieved by each agent recursively communicating its posterior probability
distribution of the target's states with neighboring agents and updating
its estimated probability distribution of the target's states using
the LogOP. As shown in Fig. \ref{fig:Demo-combine-pdfs}, combining
posterior probability distributions using the linear opinion pool
(LinOP) typically results in multimodal solutions, which are insensitive
to the weights \cite{Ref:Genest86}. On the other hand, combining
posterior probability distributions using the LogOP typically results
in unimodal, less dispersed solutions, thereby indicating a jointly
preferred consensual distribution by the network. Moreover, as discussed
in Section \ref{sub:Logarithmic-Opinion-Pool}, the optimal solution
does not depend upon the choice of scale of the prior probability
distribution and LogOP is is \textit{externally Bayesian} \cite{Ref:Genest86}
(See Fig. \ref{fig:Demo-combine-pdfs} (c-d)).

The KL divergence is the measure of the information lost when the
consensual estimate is used to approximate the locally estimated posterior
probability distributions. In Theorem \ref{thm:LogOP-minimizes-KL},
we show that the LogOP algorithm on a strongly connected (SC) balanced
graph minimizes the information lost during the consensus stage, i.e.,
the consensual probability distribution minimizes the sum of KL divergences
with the locally estimated posterior probability distributions. Methods
for communicating probability distributions and the effects of inaccuracies
on the consensual probability distribution are discussed in Section
\ref{sub:Communicating-probability-distri}.

The second contribution of this paper is the BCF algorithm presented
in Section \ref{sec:Dynamic-Bayesian-Consensus}. As illustrated in
Fig. \ref{fig:Flowchart-BCF-LogOP} and \textbf{Algorithm \ref{alg:Dynamic-Bayesian-Consensus}},
each agent generates a local estimate of the posterior probability
distribution of the target's states using the Bayesian filter. Note
that measurement exchanges with neighbors during the Bayesian filtering
stage are not mandatory and can be omitted. During the consensus stage,
the LogOP algorithm is executed multiple times to reach an agreement
across the network. The number of consensus loops ($n_{\textrm{loop}}\in\mathbb{N}$)
depends on the second largest singular value of the matrix representing
a SC balanced communication network topology. Moreover, the convergence
conditions for a given number of consensus loops are derived. Note
that this consensual probability distribution from the current time
step is used as the prior probability distribution in the next time
step, as shown in Fig \ref{fig:Flowchart-BCF-LogOP}. The novel features
of the BCF algorithm are:
\begin{itemize}
\item The algorithm can be used to track targets with general nonlinear
time-varying target dynamic models.
\item The algorithm can be used by a SC balanced network of heterogeneous
agents with general nonlinear time-varying measurement models. 
\item The algorithm achieves global exponential convergence across the network
to the consensual probability distribution of the target's states. 
\item The consensual probability distribution the best estimate, in the
information theoretic sense because it minimizes the sum of KL divergences
with the locally estimated posterior probability distributions. If
a central agent receives all the local posterior probability distributions
and is tasked to find the best estimate in the information theoretic
sense, then it would also yield the same consensual probability distribution.
Hence, we claim to have achieved distributed estimation using the
BCF algorithm. 
\end{itemize}
The Hierarchical BCF algorithms, in Section \ref{sub:Hierarchical-Bayesian-Consensus},
is used when some of the agents do not observe the target. In Section
\ref{sec:Numerical-Example}, we apply the Hierarchical BCF algorithm
to the problem of tracking orbital debris in space using the space
surveillance network on Earth.

\subsection{Notation}

The \textit{time index} is denoted by a right subscript. For example,
$\boldsymbol{x}_{k}$ represents the true states of the target at
the $k^{\textrm{th}}$ time instant. The target is always within the
compact state space $\mathcal{X}$, i.e., $\boldsymbol{x}_{k}\in\mathcal{X},\thinspace\forall k\in\mathbb{N}$.
Also, $\boldsymbol{x}_{1:k}$ represents an array of the true states
of the target from the first to the $k^{\textrm{th}}$ time instant.
The \textit{agent index} is denoted by a lower-case right superscript.
For example, $\boldsymbol{z}_{k}^{j}$ represents the measurement
taken by the $j^{\textrm{th}}$ agent at the $k^{\textrm{th}}$ time
instant. The symbol $\mathbb{P}(\cdot)$ refers to probability of
an event. $\mathcal{F}_{k}^{j}$ represents the estimated probability
density function (pdf) of the target's states over the state space
$\mathcal{X}$, by the $j^{\textrm{th}}$ agent at the $k^{\textrm{th}}$
time instant. The symbol $p(\cdot)$ also refers to pdf or probability
mass function (pmf) over the state space $\mathcal{X}$. During the
consensus stage at the $k^{\textrm{th}}$ time instant, $\mathbf{\mathcal{F}}_{k,\nu}^{j}$
represents the local pdf of the target's states by the $j^{\textrm{th}}$
agent at the $\nu^{\textrm{th}}$ consensus step and $\mathbf{\mathcal{F}}_{k}^{\star}$
represents the consensual pdf to which each $\mathbf{\mathcal{F}}_{k,\nu}^{j}$
converges. Let $\mathscr{X}$ be the Borel $\sigma$--algebra on $\mathcal{X}$. 

The communication network topology at the $k^{\textrm{th}}$ time
instant is represented by the directed time-varying graph $\mathcal{G}_{k}$,
where all the agents of the system form the set of vertices $\mathcal{V}$
(which does not change with time) and the set of directed edges is
denoted by $\mathcal{E}_{k}$. The neighbors of the $j^{\textrm{th}}$
agent at the $k^{\textrm{th}}$ time instant is the set of agents
from which the $j^{\textrm{th}}$ agent receives information at the
$k^{\textrm{th}}$ time instant and is denoted by $\mathcal{N}_{k}^{j}$,
i.e., $\ell\in\mathcal{N}_{k}^{j}$ if and only if $\overrightarrow{\ell j}\in\mathcal{E}_{k}$
for all $\ell,j\in\mathcal{V}$. The set of inclusive neighbors of
the $j^{\textrm{th}}$ agent is denoted by $\mathcal{J}_{k}^{j}:=\mathcal{N}_{k}^{j}\cup\{j\}$. 

Let $\mathbb{N}$ and $\mathbb{R}$ be the sets of natural numbers
(positive integers) and real numbers respectively. The set of all
$m$ by $n$ matrices over the field of real numbers $\mathbb{R}$
is denoted by $\mathbb{R}^{m\times n}$. Let $\lambda$ and $\sigma$
represent the eigenvalue and the singular value of a square matrix.
Let $\mathbf{1}=[1,1,\ldots,1]^{T}$, $\mathbf{I}$, $\mathbf{0}$,
and $\phi$ be the ones vector, the identity matrix, the zero matrix
of appropriate sizes, and the empty set respectively. The symbols
$\left|\cdot\right|$, $\left\lceil \cdot\right\rceil $, and $\textrm{sgn}(\cdot)$
represent the absolute value, ceiling function, and signum function
respectively. Let $\ln(\cdot)$ and $\log_{c}(\cdot)$ represent the
natural logarithm and the logarithm to the base $c$. Finally, $\|\cdot\|_{\ell_{p}}$
represents the $\ell_{p}$ vector norm. The $\mathcal{L}_{p}$ function
denotes the set of all functions $f(\boldsymbol{x}):\mathbb{R}^{n_{x}}\rightarrow\mathbb{R}$
with the bounded integral $\left(\int_{\mathcal{X}}|f(\boldsymbol{x})|^{p}d\mu(\boldsymbol{x})\right)^{1/p}$,
where $\mu$ is a measure on $\mathscr{X}$.

\section{Preliminaries \label{sec:Preliminaries}}

In this section, we first state four assumptions used throughout this
paper and then introduce the problem statement of BCF. Next, we discuss
an extension of the Bayesian filter to sensor fusion over a network. 

\begin{assumption} In this paper, all the algorithms are presented
in discrete time. \label{assump:discrete_time} \hfill $\Box$\end{assumption}

\begin{assumption} The state space ($\mathcal{X}\subset\mathbb{R}^{n_{x}}$)
is closed and bounded. Hence, by the Heine\textendash{}Borel theorem
(cf. \cite[pp. 86]{Ref:Ross80}), $\mathcal{X}$ is compact. \label{assump:compact_state_space}
\hfill $\Box$\end{assumption}

\begin{assumption} All continuous probability distributions are upper-bounded
by some large value $\mathcal{M}\in\mathbb{R}$. \label{assump:pdf-upper-bound}
\hfill $\Box$\end{assumption}

\begin{assumption} The inter-agent communication time scale is much
faster than the tracking/estimation time scale. \label{assump:comm_tracking_time_scales}
\hfill $\Box$\end{assumption}

Assumptions \ref{assump:discrete_time} and \ref{assump:compact_state_space}
are introduced to discretize the time and bound the state space, so
that the algorithms are computationally tractable. Under these assumptions,
particle filters \cite{Ref:Arulampalam02}, approximate grid--based
filters, or histogram filters \cite{Ref:Thrun05} can be used to execute
the algorithms developed in this paper. Assumptions \ref{assump:compact_state_space}
and \ref{assump:pdf-upper-bound} are introduced to take advantage
of the results in information theory and measure theory, which deal
with bounded functions on compact support. Under Assumption \ref{assump:comm_tracking_time_scales},
the agents can execute multiple consensus loops within each tracking
time step. 

We envisage that the results in this paper could be extended to continuous
time if the Fokker--Plank equations are solved efficiently \cite{Ref:Daum05}
and additional issues due to communication delay and time scale separation
are addressed. Under Assumption \ref{assump:discrete_time}, we do
not discuss continuous time related issues in this paper. Next, we
show that discrete and continuous probability distributions can be
handled in a unified manner. 

\begin{remark} Let $\mathscr{X}$ be the Borel $\sigma$--algebra
for $\mathcal{X}$. The probability of a set $\mathscr{A}\in\mathscr{X}$
may be written as the Lebesgue--Stieltjes integral $\mathbb{P}(\mathscr{A})=\int_{\mathscr{A}}p(\boldsymbol{x})\thinspace d\mu(\boldsymbol{x})$,
where $\mu$ is a measure on $\mathscr{X}$. In the continuous case,
$p(\boldsymbol{x})$ is the pdf and $\mu$ is the Lebesgue measure.
In the discrete case, $p(\boldsymbol{x})$ is the pmf and $\mu$ is
the the counting measure. \hfill $\Box$ \end{remark}

Hence, in this paper, we only deal with pdfs over $\mathcal{X}$ with
$\mu$ as the Lebesgue measure. Similar arguments will also work for
pmfs or mixed probability distributions.

\subsection{Problem Statement \label{sub:Problem-Statement}}

Let $\mathcal{X}\subset\mathbb{R}^{n_{x}}$ be the $n_{x}$-dimensional
state space of the target. The dynamics of the target in discrete
time $\{\boldsymbol{x}_{k},k\in\mathbb{N},\boldsymbol{x}_{k}\in\mathcal{X}\}$
is given by: 
\begin{equation}
\boldsymbol{x}_{k}=\boldsymbol{f}_{k}(\boldsymbol{x}_{k-1},\boldsymbol{v}_{k-1})\thinspace,\label{sys_mod}
\end{equation}
where $\boldsymbol{f}_{k}:\mathbb{R}^{n_{x}}\times\mathbb{R}^{n_{v}}\rightarrow\mathbb{R}^{n_{x}}$
is a possibly nonlinear time-varying function of the state $\boldsymbol{x}_{k-1}$
and an independent and identically distributed (i.i.d.) process noise
$\boldsymbol{v}_{k-1}$, where $n_{v}$ is the dimension of the process
noise vector. 

Let $m$ heterogeneous agents simultaneously track this target and
estimate the pdf of the target's states (where $m$ does not change
with time). The measurement model of the $j^{\textrm{th}}$ agent
is given by: 
\begin{equation}
\boldsymbol{z}_{k}^{j}=\boldsymbol{h}_{k}^{j}(\boldsymbol{x}_{k},\boldsymbol{w}_{k}^{j}),\quad\forall j\in\{1,\ldots,m\},\label{mes_mod_con}
\end{equation}
where $\boldsymbol{h}_{k}^{j}:\mathbb{R}^{n_{x}}\times\mathbb{R}^{n_{wj}}\rightarrow\mathbb{R}^{n_{zj}}$
is a possibly nonlinear time-varying function of the state $\boldsymbol{x}_{k}$
and an i.i.d. measurement noise $\boldsymbol{w}_{k}^{j}$, where $n_{zj},n_{wj}$
are dimensions of the measurement and measurement noise vectors respectively.
Note that the measurement model of agents is quite general since it
accommodates heterogeneous sensors with various bandwidths, ranges,
and noise characteristics and partial state observation. 

The objective of the BCF is to estimate the target's states and maintain
consensus across the network. This objective is achieved in two steps:
(i) each agent locally estimates the pdf of the target's states using
a Bayesian filter, and (ii) each agent's local estimate converges
to a global estimate during the consensus stage (see Fig. \ref{fig:BCF-overview}).
The objective of Bayesian filtering with/without measurement exchange,
discussed in Section \ref{sub:Bayesian-Filtering-Algorithms}, is
to estimate the posterior pdf of the target's states at the $k^{\textrm{th}}$
time instant, which is denoted by $\mathbf{\mathcal{F}}_{k}^{j},\forall j\in\{1,\ldots,m\}$,
using the estimated prior pdf of the target's states $\mathbf{\mathcal{F}}_{k-1}^{j}$
from the $(k-1)^{\textrm{th}}$ time instant and the new measurement
array obtained at the $k^{\textrm{th}}$ time instant. The objective
of the consensus stage, discussed in Section \ref{sec:Bayesian-Consensus-Filters},
is to guarantee pointwise convergence of each estimated pdf $\mathbf{\mathcal{F}}_{k,\nu}^{j}$
to the consensual pdf $\mathbf{\mathcal{F}}_{k}^{\star}$.

\subsection{Bayesian Filter with Measurement Exchange \label{sub:Bayesian-Filtering-Algorithms}}

A Bayesian filter consist of two steps: (i) the prior pdf of the target's
states is obtained during the prediction stage, and (ii) the posterior
pdf of the target's states is updated using the new measurement array
during the update stage \nocite{Ref:DeGroot75,Ref:Jeffreys61,Ref:Subrahmaniam79,Ref:Pearl88}\cite{Ref:DeGroot75}--\cite{Ref:Pearl88}.
The Bayesian filter gives the exact posterior probability distribution,
hence it is the best possible estimate of the target from the available
information. 

Exchange of measurements is optional since heterogeneous agents, with
different priors, fields of view, resolutions, tolerances, etc., may
not be able to combine measurements from other agents. For example,
if a satellite in space and a low flying quadrotor are observing the
same target, then they cannot exchange measurements due to their different
fields of view. Furthermore, a centralized estimator may not be able
to combine measurements from all heterogeneous agents in the network
to estimate $p_{k}(\boldsymbol{x}_{k}|\boldsymbol{z}_{k}^{\{1,\ldots,m\}})$,
because it would have to use a common prior for all the agents. Hence,
in this paper, we let the individual agents generate their own posterior
pdfs of the target's states and then combine them to get the best
estimated pdf from the network. 

If an agent can combine measurements from another neighboring agent
during its update stage, then we call them \textit{measurement neighbors}.
In this section, we extend the Bayesian filter by assuming that each
agent transmits its measurements to other agents in the network, and
receives the measurements from its measurement neighbors. Let $\boldsymbol{z}_{k}^{\mathcal{S}_{k}^{j}}:=\{\boldsymbol{z}_{k}^{\ell},\forall\ell\in\mathcal{S}_{k}^{j}\}$
denote the array of measurements taken at the $k^{\textrm{th}}$ time
instant by the measurement neighbors of the $j^{\textrm{th}}$ agent,
where $\mathcal{S}_{k}^{j}\subseteq\mathcal{J}_{k}^{j}$ denotes the
set of measurement neighbors among the inclusive neighbors of the
$j^{\textrm{th}}$ agent. Next, we assume that the prior is available
at the initial time. 

\begin{assumption} For each agent $j$, the initial prior of the
states $\mathbf{\mathcal{F}}_{0}^{j}=p_{0}^{j}(\boldsymbol{x}_{0})$,
is assumed to be available. In case no knowledge about $\boldsymbol{x}_{0}$
is available, $\mathbf{\mathcal{F}}_{0}^{j}$ is assumed to be uniformly
distributed over $\mathcal{X}$. \label{assump:prior} \hfill $\Box$
\end{assumption}

In\emph{ }\textit{Bayesian Filtering with Measurement Exchanges},
the $j^{\textrm{th}}$ agent estimates the posterior pdf of the target's
states $\mathbf{\mathcal{F}}_{k}^{j}=p_{k}^{j}(\boldsymbol{x}_{k}|\boldsymbol{z}_{k}^{\mathcal{S}_{k}^{j}})$
at the $k^{\textrm{th}}$ time instant using the estimated consensual
pdf of the target's states $\mathbf{\mathcal{F}}_{k-1}^{j}=p_{k-1}^{j}(\boldsymbol{x}_{k-1})$
from the $(k-1)^{\textrm{th}}$ time instant and the new measurement
array $\boldsymbol{z}_{k}^{\mathcal{S}_{k}^{j}}$ obtained at the
$k^{\textrm{th}}$ time instant. The prediction stage involves using
the target dynamics model (\ref{sys_mod}) to obtain the estimated
pdf of the target's states at the $k^{\textrm{th}}$ time instant
via the Chapman--Kolmogorov equation: \vspace{-10pt}

\begin{align}
p_{k}^{j}(\boldsymbol{x}_{k}) & =\int_{\mathcal{X}}\thinspace p_{k}^{j}(\boldsymbol{x}_{k}|\boldsymbol{x}_{k-1})\thinspace p_{k-1}^{j}(\boldsymbol{x}_{k-1})\thinspace d\mu(\boldsymbol{x}_{k-1})\label{predict_stage_con}
\end{align}
The probabilistic model of the state evolution $p_{k}^{j}(\boldsymbol{x}_{k}|\boldsymbol{x}_{k-1})$
is defined by the target dynamics model (\ref{sys_mod}) and the known
statistics of the i.i.d. process noise $\boldsymbol{v}_{k-1}$. 

\begin{proposition} The new measurement array ($\boldsymbol{z}_{k}^{\mathcal{S}_{k}^{j}}$)
is used to compute the posterior pdf of the target's states ($\mathbf{\mathcal{F}}_{k}^{j}=p_{k}^{j}(\boldsymbol{x}_{k}|\boldsymbol{z}_{k}^{\mathcal{S}_{k}^{j}})$)
during the update stage using Bayes' rule (\ref{update_stage1_con}):
\vspace{-10pt}

\begin{align}
\begin{aligned}p_{k}^{j}(\boldsymbol{x}_{k}|\boldsymbol{z}_{k}^{\mathcal{S}_{k}^{j}})=\frac{\Big(\prod_{\ell\in\mathcal{S}_{k}^{j}}p_{k}^{\ell}(\boldsymbol{z}_{k}^{\ell}|\boldsymbol{x}_{k})\Big)\thinspace p_{k}^{j}(\boldsymbol{x}_{k})}{\int_{\mathcal{X}}\Big(\prod_{\ell\in\mathcal{S}_{k}^{j}}p_{k}^{\ell}(\boldsymbol{z}_{k}^{\ell}|\boldsymbol{x}_{k})\Big)p_{k}^{j}(\boldsymbol{x}_{k})\thinspace d\mu(\boldsymbol{x}_{k})},\end{aligned}
\label{update_stage1_con}
\end{align}
The likelihood function $p_{k}^{\ell}(\boldsymbol{z}_{k}^{\ell}|\boldsymbol{x}_{k}),\thinspace\forall\ell\in\mathcal{S}_{k}^{j}$
is defined by the measurement model (\ref{mes_mod_con}), and the
corresponding known statistics of the i.i.d. measurement noise $\boldsymbol{w}_{k}^{\ell}$.
\end{proposition}
\begin{IEEEproof}
We need to show that the term $p(\boldsymbol{z}_{k}^{\mathcal{S}_{k}^{j}}|\boldsymbol{x}_{k})$
in the Bayesian filter \cite{Ref:Pearl88} simplifies to $\Big(\prod_{\ell\in\mathcal{S}_{k}^{j}}p_{k}^{\ell}(\boldsymbol{z}_{k}^{\ell}|\boldsymbol{x}_{k})\Big)$.
Let the agents $r,r+1,\ldots$ be measurement neighbors of the $j^{\textrm{th}}$
agent at the $k^{\textrm{th}}$ time instant (i.e., $r,r+1,\ldots\in\mathcal{N}_{k}^{j}$).
Let us define $\boldsymbol{z}_{k}^{\mathcal{S}_{k}^{j}\backslash\{r\}}:=\left\{ \boldsymbol{z}_{k}^{\ell},\thinspace\forall\ell\in\mathcal{S}_{k}^{j}\setminus\{r\}\right\} $
as the measurement array obtained by the $j^{\textrm{th}}$ agent
at the $k^{\textrm{th}}$ time instant, which does not contain the
measurement from the $r^{\textrm{th}}$ agent. Since the measurement
noise is i.i.d., (\ref{mes_mod_con}) describes a Markov process of
order one, we get: \vspace{-15pt}

\begin{align*}
 & \frac{p(\boldsymbol{z}_{k}^{\mathcal{S}_{k}^{j}},\boldsymbol{x}_{k})}{p(\boldsymbol{x}_{k})}\!=\!\frac{p(\boldsymbol{z}_{k}^{\mathcal{S}_{k}^{j}},\boldsymbol{x}_{k})}{p(\boldsymbol{z}_{k}^{\mathcal{S}_{k}^{j}\backslash\{r\}},\boldsymbol{x}_{k})}\frac{p(\boldsymbol{z}_{k}^{\mathcal{S}_{k}^{j}\backslash\{r\}}\!,\boldsymbol{x}_{k})}{p(\boldsymbol{z}_{k}^{\mathcal{S}_{k}^{j}\backslash\{r,r+1\}}\!,\boldsymbol{x}_{k})}\ldots\frac{p(\boldsymbol{z}_{k}^{j},\boldsymbol{x}_{k})}{p(\boldsymbol{x}_{k})}\\
 & =\! p(\boldsymbol{z}_{k}^{r}|\boldsymbol{z}_{k}^{\mathcal{S}_{k}^{j}\backslash\{r\}}\!,\boldsymbol{x}_{k})p(\boldsymbol{z}_{k}^{(r+1)}|\boldsymbol{z}_{k}^{\mathcal{S}_{k}^{j}\backslash\{r,r+1\}}\!,\boldsymbol{x}_{k})\ldots p(\boldsymbol{z}_{k}^{j}|\boldsymbol{x}_{k})
\end{align*}
Thus, we obtain $p(\boldsymbol{z}_{k}^{\mathcal{S}_{k}^{j}}|\boldsymbol{x}_{k})=p(\boldsymbol{z}_{k}^{\mathcal{S}_{k}^{j}},\boldsymbol{x}_{k})/p(\boldsymbol{x}_{k})=\prod_{\ell\in\mathcal{S}_{k}^{j}}p_{k}^{\ell}(\boldsymbol{z}_{k}^{\ell}|\boldsymbol{x}_{k})$. 
\end{IEEEproof}
\vspace{5pt} If the estimates are represented by pmfs, then we can
compare these estimates using entropy \cite[pp. 13]{Ref:Cover91},
which is a measure of the uncertainty associated with a random variable
or its information content.

\begin{remark} Let $\boldsymbol{Y}_{k}^{j}$ be a random variable
with the pmf $p_{k}^{j}(\boldsymbol{x}_{k}|\boldsymbol{z}_{k}^{j})$
given by a stand-alone Bayesian filter \cite{Ref:Pearl88}, $\boldsymbol{Y}_{k}^{\mathcal{S}_{k}^{j}}$
is a random variable with the pmf $p_{k}^{j}(\boldsymbol{x}_{k}|\boldsymbol{z}_{k}^{\mathcal{S}_{k}^{j}})$
given by Bayesian filter with measurement exchange, and $H(\cdot)$
refers to the entropy of the random variable. Since $\boldsymbol{Y}_{k}^{\mathcal{S}_{k}^{j}}$
is obtained by conditioning $\boldsymbol{Y}_{k}^{j}$ with $\boldsymbol{z}_{k}^{\mathcal{S}_{k}^{j}\backslash\{j\}}$
because $p_{k}^{j}(\boldsymbol{x}_{k}|\boldsymbol{z}_{k}^{\mathcal{S}_{k}^{j}})=p_{k}^{j}(\boldsymbol{x}_{k}|\boldsymbol{z}_{k}^{j},\boldsymbol{z}_{k}^{\mathcal{S}_{k}^{j}\backslash\{j\}})$,
the claim follows from the theorem on conditioning reduces entropy
(cf. \cite[pp. 27]{Ref:Cover91}): 
\begin{align*}
0\leq I(\boldsymbol{Y}_{k}^{j};\boldsymbol{z}_{k}^{\mathcal{S}_{k}^{j}\backslash\{j\}})=H(\boldsymbol{Y}_{k}^{j})-H(\boldsymbol{Y}_{k}^{j}|\boldsymbol{z}_{k}^{\mathcal{S}_{k}^{j}\backslash\{j\}})\thinspace,
\end{align*}
where $I(\cdot;\cdot)$ refers to the nonnegative mutual information
between two random variables. Since $H(\boldsymbol{Y}_{k}^{\mathcal{S}_{k}^{j}})=H(\boldsymbol{Y}_{k}^{j}|\boldsymbol{z}_{k}^{\mathcal{S}_{k}^{j}\backslash\{j\}})\leq H(\boldsymbol{Y}_{k}^{j})$,
the estimate given by Bayesian filter with measurement exchange is
more accurate than that obtained by a stand-alone Bayesian filter.
\hfill $\Box$ \end{remark}

Note that (\ref{update_stage1_con}) is similar to the empirical equation
for \textit{Independent Likelihood Pool} given in \cite{Ref:Punska99}
and a generalization of the \textit{Distributed Sequential Bayesian
Estimation Algorithm} given in \cite{Ref:Nehorai07}. The structure
of (\ref{update_stage1_con}) ensures that an arbitrary part of the
prior distribution does not dominate the measurements. There is no
consensus protocol across the network because each agent receives
information only from its neighboring agents and never receives measurements
(even indirectly) from any other agent in the network.

\section{Combining Probability Distributions \label{sec:Bayesian-Consensus-Filters}}

In this section, we present the algorithms for achieving consensus
in probability distributions across the network. As discussed before,
the objective of the consensus stage in \textbf{Algorithm \ref{alg:Dynamic-Bayesian-Consensus}}
is to guarantee pointwise convergence of each $\mathbf{\mathcal{F}}_{k}^{j}$
to a consensual pdf $\mathbf{\mathcal{F}}_{k}^{\star}$, which is
independent of $j$. This is achieved by each agent recursively transmitting
its estimated pdf of the target's states to other agents, receiving
estimates of its neighboring agents, and updating its estimate of
the target. Let $\mathbf{\mathcal{F}}_{k,0}^{j}=\mathbf{\mathcal{F}}_{k}^{j}$
represent the local estimated posterior pdf of the target's states,
by the $j^{\textrm{th}}$ agent at the start of the consensus stage,
obtained using Bayesian filters with/without measurement exchange.
During each of the $n_{\textrm{loop}}$ iterations within the consensus
stage in \textbf{Algorithm \ref{alg:Dynamic-Bayesian-Consensus}},
this estimate is updated as follows: \vspace{-10pt}

\begin{equation}
\mathbf{\mathcal{F}}_{k,\nu}^{j}\!=\!\mathcal{T}\left(\cup_{\ell\in\mathcal{J}_{k}^{j}}\{\mathbf{\mathcal{F}}_{k,\nu-1}^{\ell}\}\right),\forall j\in\{1,\ldots,m\},\forall\nu\in\mathbb{N},\label{combining_estimates}
\end{equation}
where $\mathcal{T}(\cdot)$ is the linear or logarithmic opinion pool
for combining the pdf estimates. Note that the problem of measurement
neighbors does not arise here since all pdfs are expressed over the
complete state space $\mathcal{X}$. 

We introduce Lemma \ref{lemma:convegence-distribution-measure} to
show that pointwise convergence of pdfs is the sufficient condition
for convergence of their induced measures in total variation (TV).
Let $\mathcal{F}_{1},\ldots,\lim_{n\rightarrow\infty}\mathcal{F}_{n},\thinspace\mathcal{F}^{\star}$
be real-valued measurable functions on $\mathcal{X}$, $\mathscr{X}$
be the Borel $\sigma$-algebra of $\mathcal{X}$, and $\mathscr{A}$
be any set in $\mathscr{X}$. If $\mu_{\mathcal{F}_{n}}(\mathscr{A})=\int_{\mathscr{A}}\mathcal{F}_{n}(\boldsymbol{x})d\mu(\boldsymbol{x})$
for any set $\mathscr{A}\in\mathscr{X}$, then $\mu_{\mathcal{F}_{n}}$
is defined as the measure induced by the function $\mathcal{F}_{n}$
on $\mathscr{X}$. Let $\mu_{\mathcal{F}_{n}},\thinspace\mu_{\mathcal{F}^{\star}}$
denote the respective induced measures of $\mathcal{F}_{n},\thinspace\mathcal{F}^{\star}$
on $\mathscr{X}$.

\begin{definition} \textit{(Convergence in TV)} If $\|\mu_{\mathcal{F}_{n}}-\mu_{\mathcal{F}^{\star}}\|_{\textrm{TV}}:=\sup_{\mathscr{A}\in\mathscr{X}}|\mu_{\mathcal{F}_{n}}(\mathscr{A})-\mu_{\mathcal{F}^{\star}}(\mathscr{A})|$
tends to zero as $n\rightarrow\infty$, then the measure $\mu_{\mathcal{F}_{n}}$
converges to the measure $\mu_{\mathcal{F}^{\star}}$ in TV, i.e.,
$\lim_{n\rightarrow\infty}\mu_{\mathcal{F}_{n}}\xrightarrow{\textrm{T.V.}}\mu_{\mathcal{F}^{\star}}$.
\hfill $\Box$ \end{definition}

\begin{lemma} \label{lemma:convegence-distribution-measure} \textit{(Pointwise
convergence implies convergence in TV)} If $\mathcal{F}_{n}$ converges
to $\mathcal{F}^{\star}$ pointwise, i.e., $\lim_{n\rightarrow\infty}\mathcal{F}_{n}=\mathcal{F}^{\star}$
pointwise; then the measure $\mu_{\mathcal{F}_{n}}$ converges in
TV to the measure $\mu_{\mathcal{F}^{\star}}$, i.e., $\lim_{n\rightarrow\infty}\mu_{\mathcal{F}_{n}}\xrightarrow{\textrm{T.V.}}\mu_{\mathcal{F}^{\star}}$.
\end{lemma}
\begin{IEEEproof}
Similar to the proof of Scheff$\acute{\textrm{e}}$'s theorem \cite[pp. 84]{Ref:Durrett05},
under Assumption \ref{assump:pdf-upper-bound}, using the dominated
convergence theorem (cf. \cite[Theorem 1.5.6, pp. 23]{Ref:Durrett05})
for any set $\mathscr{A}\in\mathscr{X}$ gives: \vspace{-10pt}

\begin{align*}
\lim_{n\rightarrow\infty}\!\int_{\mathscr{A}}\!\!\!\mathcal{F}_{n}(\boldsymbol{x})d\mu(\boldsymbol{x}) & \!=\!\int_{\mathscr{A}}\!\lim_{n\rightarrow\infty}\!\mathcal{F}_{n}(\boldsymbol{x})d\mu(\boldsymbol{x})\!=\!\int_{\mathscr{A}}\!\!\!\mathcal{F}^{\star}(\boldsymbol{x})d\mu(\boldsymbol{x}).
\end{align*}
This relation between measures implies that $\|\lim_{n\rightarrow\infty}\mu_{\mathcal{F}_{n}}-\mu_{\mathcal{F}^{\star}}\|_{\textrm{TV}}=0$
and $\lim_{n\rightarrow\infty}\mu_{\mathcal{F}_{n}}\xrightarrow{\textrm{T.V.}}\mu_{\mathcal{F}^{\star}}$.
\end{IEEEproof}

\subsection{Consensus using the Linear Opinion Pool \label{sub:Linear-Opinion-Pool}}

The first method of combining the estimates is motivated by the linear
consensus algorithms widely studied in the literature \nocite{Ref:Ren05TAC,Ref:Saber07,Ref:Tsitsiklis05}\cite{Ref:Ren05TAC}--\cite{Ref:Tsitsiklis05}.
The pdfs are combined using the Linear Opinion Pool (LinOP) of probability
measures \cite{Ref:DeGroot74,Ref:Genest86}: \vspace{-10pt} 

\begin{align}
\mathbf{\mathcal{F}}_{k,\nu}^{j}= & \sum_{\ell\in\mathcal{J}_{k}^{j}}a_{k,\nu-1}^{j\ell}\mathbf{\mathcal{F}}_{k,\nu-1}^{\ell},\thinspace\forall j\in\{1,\ldots,m\},\forall\nu\in\mathbb{N},\label{agreement_equation}
\end{align}
where $\sum_{\ell\in\mathcal{J}_{k}^{j}}a_{k,\nu-1}^{j\ell}=1$ and
the updated pdf $\mathbf{\mathcal{F}}_{k,\nu}^{j}$ after the $\nu^{\textrm{th}}$
consensus loop is a weighted average of the pdfs of the inclusive
neighbors $\mathbf{\mathcal{F}}_{k,\nu-1}^{\ell},\forall\ell\in\mathcal{J}_{k}^{j}$
from the $(\nu-1)^{\textrm{th}}$ consensus loop, at the $k^{\textrm{th}}$
time instant. Let $\mathbf{\mathcal{W}}_{k,\nu}:=\left(\mathbf{\mathcal{F}}_{k,\nu}^{1},\ldots,\mathbf{\mathcal{F}}_{k,\nu}^{m}\right)^{T}$
denote an array of pdf estimates of all the agents after the $\nu^{\textrm{th}}$
consensus loop, then the LinOP (\ref{agreement_equation}) can be
expressed concisely as: 
\begin{equation}
\mathbf{\mathcal{W}}_{k,\nu}=P_{k,\nu-1}\mathbf{\mathcal{W}}_{k,\nu-1},\thinspace\forall\nu\in\mathbb{N},\label{PW_eqn}
\end{equation}
where $P_{k,\nu-1}$ is a matrix with entries $P_{k,\nu-1}[j,\ell]=a_{k,\nu-1}^{j\ell}$.

\begin{assumption} The communication network topology of the multi-agent
system $\mathcal{G}_{k}$ is strongly connected (SC). The weights
$a_{k,\nu-1}^{j\ell},\forall j,\ell\in\{1,\ldots,m\}$ and the matrix
$P_{k,\nu-1}$ have the following properties: (i) the weights are
the same for all consensus loops within each time instant, i.e., $a_{k,\nu-1}^{j\ell}=a_{k}^{j\ell}$
and $P_{k,\nu-1}=P_{k},\forall\nu\in\mathbb{N}$; (ii) the matrix
$P_{k}$ conforms with the graph $\mathcal{G}_{k}$, i.e., $a_{k}^{j\ell}>0$
if and only if $\ell\in\mathcal{J}_{k}^{j}$, else $a_{k}^{j\ell}=0$;
and (iii) the matrix $P_{k}$ is row stochastic, i.e., $\sum_{\ell=1}^{m}a_{k}^{j\ell}=1$.
\label{assump:weights1} \hfill $\Box$ \end{assumption}

\begin{theorem} \label{thm:LinOP-SC-digraphs} \emph{(Consensus using
the LinOP on SC Digraphs)} Under Assumption \ref{assump:weights1},
using the LinOP (\ref{agreement_equation}), each $\mathbf{\mathcal{F}}_{k,\nu}^{j}$
asymptotically converges pointwise to the pdf $\mathbf{\mathcal{F}}_{k}^{\star}=\sum_{i=1}^{m}\pi_{i}\mathbf{\mathcal{F}}_{k,0}^{i}$
where $\boldsymbol{\pi}=[\pi_{1},\ldots\pi_{m}]^{T}$ is the unique
stationary distribution of $P_{k}$. Furthermore, their induced measures
converge in total variation, i.e., $\lim_{\nu\rightarrow\infty}\mu_{\mathbf{\mathcal{F}}_{k,\nu}^{j}}\xrightarrow{\textrm{T.V.}}\mu_{\mathcal{F}_{k}^{\star}},\thinspace\forall j\in\{1,\ldots,m\}$.
\end{theorem}
\begin{IEEEproof}
See Appendix \ref{sec:Proof-of-Theorem-BCF-LinOP-Static}.
\end{IEEEproof}
Theorem \ref{thm:LinOP-SC-digraphs} is a generalization of the linear
consensus algorithm for combining joint measurement probabilities
\cite{Ref:Rus12}. Moreover, if $\boldsymbol{\pi}=\boldsymbol{1}$
and each $\mathbf{\mathcal{F}}_{k,0}^{j}$ is a $\mathcal{L}_{2}$
function, then $\mathbf{\mathcal{F}}_{k}^{\star}=\frac{1}{m}\sum_{i=1}^{m}\mathbf{\mathcal{F}}_{k,0}^{i}$
globally minimizes the sum of the squares of $\mathcal{L}_{2}$ distances
with the locally estimated posterior pdfs. 

\begin{figure}
\begin{centering}
\begin{tabular}{cc}
\includegraphics[bb=0bp 0bp 375bp 303bp,clip,width=1.5in]{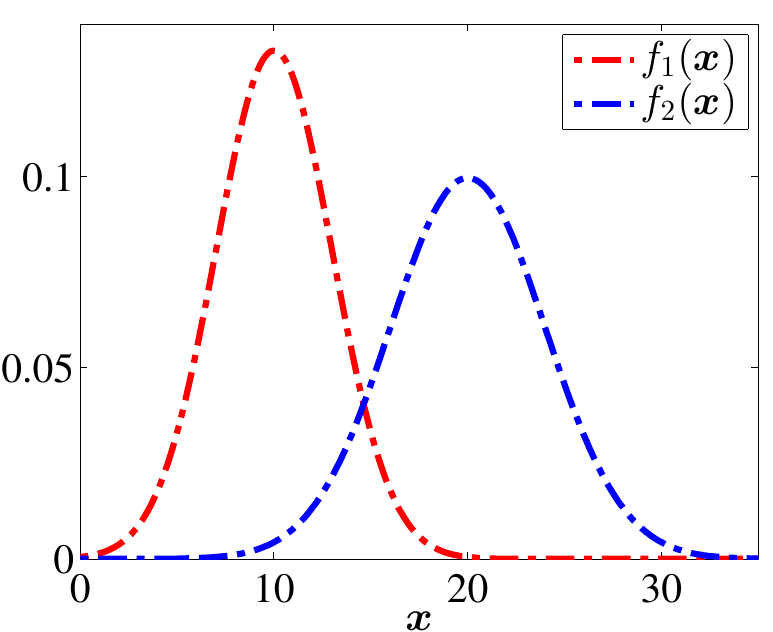} & \includegraphics[bb=0bp 0bp 375bp 303bp,clip,width=1.5in]{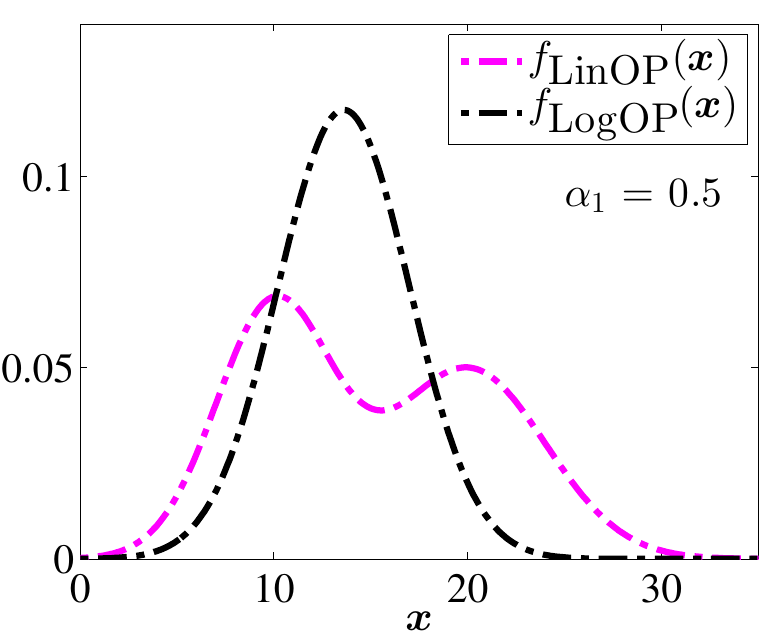}\tabularnewline
(a) & (b)\tabularnewline
\includegraphics[bb=0bp 0bp 397bp 303bp,clip,width=1.6in]{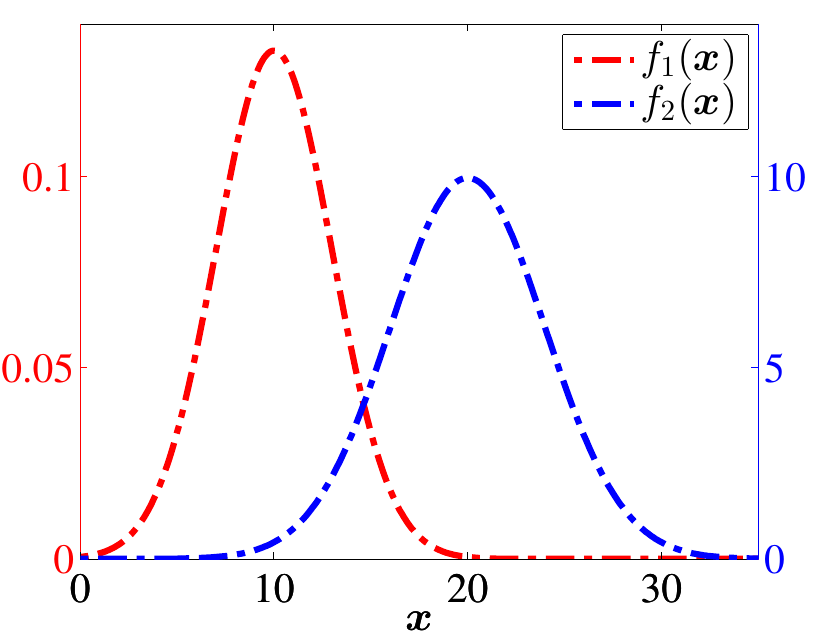} & \includegraphics[bb=0bp 0bp 371bp 303bp,clip,width=1.5in]{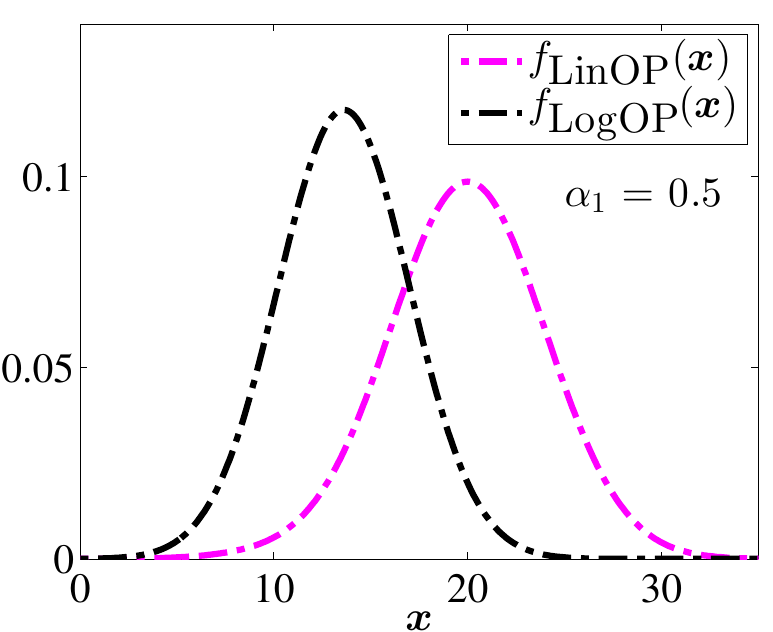}\tabularnewline
(c) & (d)\tabularnewline
\includegraphics[bb=0bp 0bp 383bp 306bp,clip,width=1.5in]{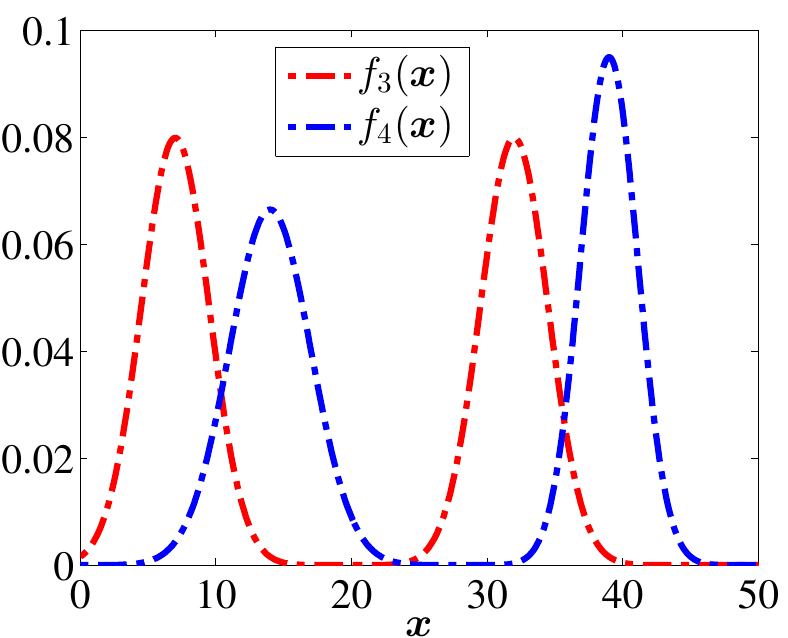} & \includegraphics[bb=0bp 0bp 378bp 306bp,clip,width=1.5in]{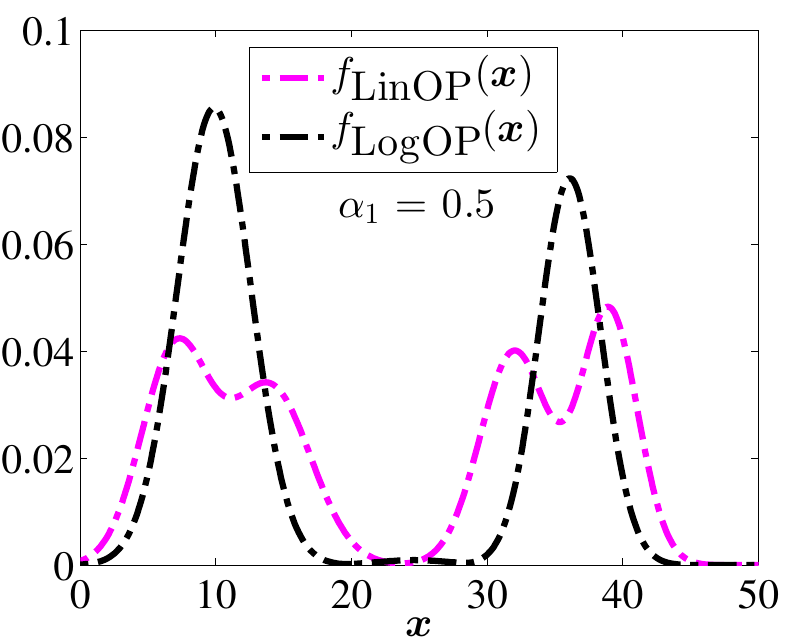}\tabularnewline
(e) & (f)\tabularnewline
\end{tabular}
\par\end{centering}

\caption{In (a), two unimodal pdfs $f_{1}(\boldsymbol{x})$ and $f_{2}(\boldsymbol{x})$
are shown. In (b), these pdfs are combined using the LinOP and LogOP
using the weight $\alpha_{1}=0.5$, i.e., $f_{\textrm{LinOP}}(\boldsymbol{x})=\left(\alpha_{1}f_{1}(\boldsymbol{x})+(1-\alpha_{1})f_{2}(\boldsymbol{x})\right)$
and $f_{\textrm{LogOP}}(\boldsymbol{x})=\left(f_{1}^{\alpha_{1}}\times f_{2}^{(1-\alpha_{1})}\right)/\left(\int_{\mathcal{X}}\left(f_{1}^{\alpha_{1}}\times f_{2}^{(1-\alpha_{1})}\right)d\mu(\boldsymbol{x})\right)$.
Note that the LinOP solution is multimodal while the LogOP solution
is unimodal, indicating a consensual pdf. \hfill\null \hspace{200pt}
In (c), the scale of the function $f_{2}(\boldsymbol{x})$ is changed
to $0$-$100$ from the standard $0$-$1$ scale. In (d), the normalized
LinOP solution changes drastically but the LogOP solution remains
unaffected. \hfill\null \hspace{150pt} In (e), the pdfs $f_{3}(\boldsymbol{x})$
and $f_{4}(\boldsymbol{x})$ have bimodal nature. In (f), the LogOP
solution preserves this bimodal nature. \label{fig:Demo-combine-pdfs}
\vspace{-10pt}}
\end{figure}

As shown in Fig. \ref{fig:Demo-combine-pdfs} (a-b), the main difficulty
with the LinOP is that the resulting solution is typically multimodal,
so no clear choice for jointly preferred estimate emerges from it
\cite{Ref:Genest86}. Moreover, the LinOP algorithm critically depends
on the assumption that the same $0$-$1$ scale is used by every agent
as shown in Fig. \ref{fig:Demo-combine-pdfs} (c-d). Hence, better
schemes for combining probability distributions are needed for the
proposed BCF algorithm. 

\subsection{Consensus using the Logarithmic Opinion Pool \label{sub:Logarithmic-Opinion-Pool}}

Note that $\mathbf{\mathcal{F}}_{k,\nu}^{j}=p_{k,\nu}^{j}(\boldsymbol{x}_{k}),\forall\boldsymbol{x}_{k}\in\mathcal{X}$
represents the pdf of the estimated target's states by the $j^{\textrm{th}}$
agent during the $\nu^{\textrm{th}}$ consensus loop at the $k^{\textrm{th}}$
time instant. The LogOP is given as \cite{Ref:Bacharach79}: \vspace{-15pt}

\begin{align}
\mathcal{\mathbf{\mathcal{F}}}_{k,\nu}^{j}= & p_{k,\nu}^{j}(\boldsymbol{x}_{k})=\frac{\prod_{\ell\in\mathcal{J}_{k}^{j}}\left(p_{k,\nu-1}^{\ell}(\boldsymbol{x}_{k})\right)^{a_{k,\nu-1}^{j\ell}}}{\int_{\mathcal{X}}\prod_{\ell\in\mathcal{J}_{k}^{j}}\left(p_{k,\nu-1}^{\ell}(\boldsymbol{x}_{k})\right)^{a_{k,\nu-1}^{j\ell}}\: d\mu(\boldsymbol{x}_{k})},\nonumber \\
 & \forall j\in\{1,\ldots,m\},\forall\nu\in\mathbb{N},\label{log_agreement_equation}
\end{align}
where $\sum_{\ell\in\mathcal{J}_{k,\nu-1}^{j}}a_{k,\nu-1}^{j\ell}=1$
and the integral in the denominator of (\ref{log_agreement_equation})
is finite. Thus the updated pdf $\mathbf{\mathcal{F}}_{k,\nu}^{j}$
after the $\nu^{\textrm{th}}$ consensus loop is the weighted geometric
average of the pdfs of the inclusive neighbors $\mathbf{\mathcal{F}}_{k,\nu-1}^{\ell},\forall\ell\in\mathcal{J}_{k}^{j}$
from the $(\nu-1)^{\textrm{th}}$ consensus loop, at the $k^{\textrm{th}}$
time instant. As shown in Fig. \ref{fig:Demo-combine-pdfs} (a-b),
the LogOP solution is typically unimodal and less dispersed, indicating
a consensual estimate jointly preferred by the network \cite{Ref:Genest86}.
As shown in Fig. \ref{fig:Demo-combine-pdfs} (c-d), the LogOP solution
is invariant under under rescaling of individual degrees of belief,
hence it preserves an important credo of uni--Bayesian decision theory;
i.e., the optimal decision should not depend upon the choice of scale
for the utility function or prior probability distribution \cite{Ref:Weerahandi83}.
When the parameter space is finite and a $0$-$1$ probability scale
is adopted, the LogOP is equivalent to the Nash product \cite{Ref:Nash50}.
Note that if the local probability distribution of the target's states
is inherently multimodal, as shown in Fig. \ref{fig:Demo-combine-pdfs}
(e-f), then LogOP preserves this multimodal nature while combining
these local estimates. The most compelling reason for using LogOP
is that it is \textit{externally Bayesian}; i.e., finding the consensus
distribution commutes with the process of revising distributions using
a commonly agreed likelihood distribution. Thus \cite{Ref:Genest86}:
\vspace{-10pt}

\begin{align}
 & \mathcal{T}\left(\cup_{\ell\in\mathcal{J}_{k}^{j}}\left\{ \frac{l(\boldsymbol{x})p_{k}^{\ell}(\boldsymbol{x})}{\int_{\mathcal{X}}l(\boldsymbol{x})p_{k}^{\ell}(\boldsymbol{x})d\mu(\boldsymbol{x})}\right\} \right)\nonumber \\
 & \quad=\frac{l(\boldsymbol{x})\mathcal{T}\left(\cup_{\ell\in\mathcal{J}_{k}^{j}}\{p_{k}^{\ell}(\boldsymbol{x})\}\right)}{\int_{\mathcal{X}}l(\boldsymbol{x})\mathcal{T}\left(\cup_{\ell\in\mathcal{J}_{k}^{j}}\{p_{k}^{\ell}(\boldsymbol{x})\}\right)d\mu(\boldsymbol{x})}\thinspace,
\end{align}
where $\mathcal{T}(\cdot)$ refers to the LogOP (\ref{log_agreement_equation}),
$p_{k}^{\ell}(\boldsymbol{x}),\forall\ell\in\mathcal{J}_{k}^{j}$
are pdfs on $\mathcal{X}$ and $l(\boldsymbol{x})$ is an arbitrary
likelihood pdf on $\mathcal{X}$. Due to these advantages, LogOP is
used for combining prior distributions \cite{Ref:Rufo12} and conditional
random fields for natural language processing tasks \cite{Ref:Smith05}.
Next, we present consensus theorems using the LogOP. 

\begin{assumption} The local estimated pdf at the start of the consensus
stage is positive everywhere, i.e., $\mathcal{F}_{k,0}^{j}=p_{k,0}^{j}(\boldsymbol{x}_{k})>0,\thinspace\forall\boldsymbol{x}_{k}\in\mathcal{X},\forall j\in\{1,\ldots,m\}$.
\label{assump:nonnegative_pdf} \hfill $\Box$ \end{assumption}

Assumption \ref{assump:nonnegative_pdf} is introduced to avoid regions
with zero probability, since they would constitute vetoes and unduly
great emphasis would get placed on them. Moreover, the LogOP guarantees
that $\mathcal{F}_{k,\nu}^{j}$ will remain positive for all subsequent
consensus loop. 

\begin{definition} \textit{($\mathcal{H}_{k,\nu}^{j}$ vector for
LogOP)} For the purpose of analysis, let us choose $\boldsymbol{\psi}_{k}\in\mathcal{X}$
such that $p_{k,\nu}^{j}(\boldsymbol{\psi}_{k})>0,\forall j\in\{1,\ldots,m\},\forall\nu\in\mathbb{N}$.
Let us define $\mathbf{\mathcal{H}}_{k,\nu}^{j}:=\ln\left[\frac{p_{k,\nu}^{j}(\boldsymbol{x}_{k})}{p_{k,\nu}^{j}(\boldsymbol{\psi}_{k})}\right]$.
Under Assumption \ref{assump:nonnegative_pdf}, $\mathbf{\mathcal{H}}_{k,\nu}^{j}$
is a well-defined function, but need not be a $\mathcal{L}_{1}$ function.
Then, by simple algebraic manipulation of (\ref{log_agreement_equation}),
we get \cite{Ref:Gilardoni93}: \vspace{-10pt}

\begin{align}
 & \frac{p_{k,\nu}^{j}(\boldsymbol{x}_{k})}{p_{k,\nu}^{j}(\boldsymbol{\psi}_{k})}=\frac{\left(\frac{\prod_{\ell\in\mathcal{J}_{k}^{j}}\left(p_{k,\nu-1}^{\ell}(\boldsymbol{x}_{k})\right)^{a_{k,\nu-1}^{j\ell}}}{\int_{\mathcal{X}}\prod_{\ell\in\mathcal{J}_{k}^{j}}\left(p_{k,\nu-1}^{\ell}(\boldsymbol{x}_{k})\right)^{a_{k,\nu-1}^{j\ell}}\: d\mu(\boldsymbol{x}_{k})}\right)}{\left(\frac{\prod_{\ell\in\mathcal{J}_{k}^{j}}\left(p_{k,\nu-1}^{\ell}(\boldsymbol{\psi}_{k})\right)^{a_{k,\nu-1}^{j\ell}}}{\int_{\mathcal{X}}\prod_{\ell\in\mathcal{J}_{k}^{j}}\left(p_{k,\nu-1}^{\ell}(\boldsymbol{x}_{k})\right)^{a_{k,\nu-1}^{j\ell}}\: d\mu(\boldsymbol{x}_{k})}\right)}\thinspace,\nonumber \\
 & \mathbf{\mathcal{H}}_{k,\nu}^{j}=\sum_{\ell\in\mathcal{J}_{k}^{j}}a_{k,\nu-1}^{j\ell}\mathbf{\mathcal{H}}_{k,\nu-1}^{\ell},\thinspace\forall j\in\{1,\ldots,m\},\nu\in\mathbb{N}.\label{log2_PW_eqn}
\end{align}
Note that (\ref{log2_PW_eqn}) is similar to the LinOP (\ref{agreement_equation}).
Let $\mathbf{\mathcal{U}}_{k,\nu}:=\left(\mathbf{\mathcal{H}}_{k,\nu}^{1},\ldots,\mathbf{\mathcal{H}}_{k,\nu}^{m}\right)^{T}$
be an array of the estimates of all the agents during the $\nu^{\textrm{th}}$
consensus loop at the $k^{\textrm{th}}$ time instant, then the equation
(\ref{log2_PW_eqn}) can be expressed concisely as: 
\begin{equation}
\mathbf{\mathcal{U}}_{k,\nu}=P_{k,\nu-1}\mathbf{\mathcal{U}}_{k,\nu-1},\thinspace\forall\nu\in\mathbb{N},\label{log2_PW_eqn2}
\end{equation}
where $P_{k,\nu-1}$ is a matrix with entries $a_{k,\nu-1}^{jl}$.
\hfill $\Box$ \end{definition}

Thus we are able to use the highly nonlinear LogOP for combining the
pdf estimates, but we have reduced the complexity of the problem to
that of consensus using the LinOP. 

\begin{theorem} \label{thm:LogOP-SC-digraphs} \emph{(Consensus using
the LogOP on SC Digraphs)} Under Assumptions \ref{assump:weights1}
and \ref{assump:nonnegative_pdf}, using the LogOP (\ref{log_agreement_equation}),
each $\mathbf{\mathcal{F}}_{k,\nu}^{j}$ asymptotically converges
pointwise to the pdf $\mathbf{\mathcal{F}}_{k}^{\star}$ given by:
\vspace{-10pt}

\begin{equation}
\mathbf{\mathcal{F}}_{k}^{\star}=p_{k}^{\star}(\boldsymbol{x}_{k})=\frac{\prod_{i=1}^{m}\left(p_{k,0}^{i}(\boldsymbol{x}_{k})\right)^{\pi_{i}}}{\int_{\mathcal{X}}\prod_{i=1}^{m}\left(p_{k,0}^{i}(\boldsymbol{x}_{k})\right)^{\pi_{i}}\: d\mu(\boldsymbol{x}_{k})},\label{log2_PW_eqn4}
\end{equation}
where $\boldsymbol{\pi}$ is the unique stationary distribution of
$P_{k}$. Furthermore, their induced measures converge in total variation,
i.e., $\lim_{\nu\rightarrow\infty}\mu_{\mathbf{\mathcal{F}}_{k,\nu}^{j}}\xrightarrow{\textrm{T.V.}}\mu_{\mathcal{F}_{k}^{\star}},\thinspace\forall j\in\{1,\ldots,m\}$.
\end{theorem}
\begin{IEEEproof}
Similar to the proof of Theorem \ref{thm:LinOP-SC-digraphs}, each
$\mathbf{\mathcal{H}}_{k,\nu}^{j}$ converges pointwise to $\mathbf{\mathcal{H}}_{k}^{\star}=\boldsymbol{\pi}^{T}\mathbf{\mathcal{U}}_{k,0}=\sum_{i=1}^{m}\pi_{i}\mathbf{\mathcal{H}}_{k,0}^{i}$
asymptotically. We additionally need to show that convergence of $\mathbf{\mathcal{H}}_{k,\nu}^{j}$
to $\mathbf{\mathcal{H}}_{k}^{\star}$ implies pointwise convergence
of $\mathbf{\mathcal{F}}_{k,\nu}^{j}$ to $\mathbf{\mathcal{F}}_{k}^{\star}$.
We have $\forall\boldsymbol{x}_{k}\in\mathcal{X}$: \vspace{-10pt}

\begin{equation}
\lim_{\nu\rightarrow\infty}\!\left(\ln p_{k,\nu}^{j}(\boldsymbol{x}_{k})\!-\!\ln p_{k,\nu}^{j}(\boldsymbol{\psi}_{k})\right)\!=\!\ln p_{k}^{\star}(\boldsymbol{x}_{k})\!-\!\ln p_{k}^{\star}(\boldsymbol{\psi}_{k}).\label{simplified_H_to_F_0}
\end{equation}
We claim $\exists\bar{\boldsymbol{\psi}}_{k}\in\mathcal{X}$ such
that $\lim_{\nu\rightarrow\infty}p_{k,\nu}^{j}(\bar{\boldsymbol{\psi}}_{k})=p_{k}^{\star}(\bar{\boldsymbol{\psi}}_{k})$.
If this claim is untrue, then $0<\lim_{\nu\rightarrow\infty}p_{k,\nu}^{j}(\boldsymbol{x}_{k})<p_{k}^{\star}(\boldsymbol{x}_{k}),\forall\boldsymbol{x}_{k}\in\mathcal{X}$
or vice versa. Hence $\int_{\mathcal{X}}\lim_{\nu\rightarrow\infty}p_{k,\nu}^{j}(\boldsymbol{x}_{k})d\mu(\boldsymbol{x}_{k})=1<\int_{\mathcal{X}}p_{k}^{\star}(\boldsymbol{x}_{k})d\mu(\boldsymbol{x}_{k})$,
which results in contradiction since $p_{k}^{\star}(\boldsymbol{x}_{k})$
is also a pdf. Hence, substituting $\bar{\boldsymbol{\psi}}_{k}$
into equation (\ref{simplified_H_to_F_0}) gives $\lim_{\nu\rightarrow\infty}p_{k,\nu}^{j}(\boldsymbol{x}_{k})=p_{k}^{\star}(\boldsymbol{x}_{k}),\forall\boldsymbol{x}_{k}\in\mathcal{X}$.
Thus each $\mathbf{\mathcal{F}}_{k,\nu}^{j}$ converges pointwise
to the consensual pdf $\mathbf{\mathcal{F}}_{k}^{\star}$ given by
(\ref{log2_PW_eqn4}). By Lemma \ref{lemma:convegence-distribution-measure},
the measure induced by $\mathbf{\mathcal{F}}_{k,\nu}^{j}$ on $\mathscr{X}$
converges in total variation to the measure induced by $\mathbf{\mathcal{F}}_{k}^{\star}$
on $\mathscr{X}$, i.e., $\lim_{\nu\rightarrow\infty}\mu_{\mathbf{\mathcal{F}}_{k,\nu}^{j}}\xrightarrow{\textrm{T.V.}}\mu_{\mathbf{\mathcal{F}}_{k}^{\star}}$. 
\end{IEEEproof}
Since Perron--Frobenius theorem only yields asymptotic convergence,
we next discuss the algorithm for achieving global exponential convergence
using balanced graphs. 

\begin{assumption} In addition to Assumption \ref{assump:weights1},
the weights $a_{k}^{j\ell}$ are such that the digraph $\mathcal{G}_{k}$
is balanced. Hence for every vertex, the in-degree equals the out-degree,
i.e., $\sum_{\ell\in\mathcal{J}_{k}^{j}}a_{k}^{j\ell}=\sum_{r|j\in\mathcal{J}_{k}^{r}}a_{k}^{rj}$,
where $j,\ell,r\in\{1,\ldots,m\}$. \label{assump:weights2} \hfill
$\Box$ \end{assumption}

\begin{theorem} \label{thm:LogOP-SC-balanced-digraphs} \emph{(Consensus
using the LogOP on SC Balanced Digraphs)} Under Assumption \ref{assump:nonnegative_pdf}
and \ref{assump:weights2}, using the LogOP (\ref{log_agreement_equation}),
each $\mathbf{\mathcal{F}}_{k,\nu}^{j}$ globally exponentially converges
pointwise to the pdf $\mathbf{\mathcal{F}}_{k}^{\star}$ given by:
\vspace{-15pt}

\begin{equation}
\mathbf{\mathcal{F}}_{k}^{\star}=p_{k}^{\star}(\boldsymbol{x}_{k})=\frac{\prod_{i=1}^{m}\left(p_{k,0}^{i}(\boldsymbol{x}_{k})\right)^{\frac{1}{m}}}{\int_{\mathcal{X}}\prod_{i=1}^{m}\left(p_{k,0}^{i}(\boldsymbol{x}_{k})\right)^{\frac{1}{m}}\: d\mu(\boldsymbol{x}_{k})}\label{log2_PW_eqn5}
\end{equation}
at a rate faster or equal to $\sqrt{\lambda_{m-1}(P_{k}^{T}P_{k})}=\sigma_{m-1}(P_{k})$.
Furthermore, their induced measures globally exponentially converge
in total variation, i.e., $\lim_{\nu\rightarrow\infty}\mu_{\mathbf{\mathcal{F}}_{k,\nu}^{j}}\xrightarrow{\textrm{T.V.}}\mu_{\mathcal{F}_{k}^{\star}},\thinspace\forall j\in\{1,\ldots,m\}$.
\end{theorem}
\begin{IEEEproof}
Since Assumption \ref{assump:weights2} is stronger than Assumption
\ref{assump:weights1}, we get $\lim_{\nu\rightarrow\infty}P_{k}^{\nu}=\mathbf{1}\boldsymbol{\pi}^{T}$.
Moreover, since $P_{k}$ is also a column stochastic matrix, therefore
$\boldsymbol{\pi}=\frac{1}{m}\mathbf{1}$ is its left eigenvector
corresponding to the eigenvalue $1$, i.e., $P_{k}^{T}\frac{1}{m}\mathbf{1}=1\frac{1}{m}\mathbf{1}$
and satisfying the normalizing condition. Hence, we get $\lim_{\nu\rightarrow\infty}P_{k}^{\nu}=\frac{1}{m}\mathbf{1}\mathbf{1}^{T}$
and each $\mathbf{\mathcal{H}}_{k,\nu}^{j}$ converges pointwise to
$\mathbf{\mathcal{H}}_{k}^{\star}=\frac{1}{m}\mathbf{1}^{T}\mathbf{\mathcal{U}}_{k,0}=\frac{1}{m}\sum_{i=1}^{m}\mathbf{\mathcal{H}}_{k,0}^{i}$. 

From the proof of Theorem \ref{thm:LogOP-SC-digraphs}, we get that
each $\mathbf{\mathcal{F}}_{k,\nu}^{j}$ converges pointwise to the
consensual pdf $\mathbf{\mathcal{F}}_{k}^{\star}$ given by (\ref{log2_PW_eqn5}).
Note that $\mathbf{\mathcal{F}}_{k,\nu}^{j}$ are $\mathcal{L}_{1}$
functions but $\mathbf{\mathcal{H}}_{k,\nu}^{j}$ need not be $\mathcal{L}_{1}$
functions. Let $V_{\textrm{tr}}=\left[\frac{1}{\sqrt{m}}\mathbf{1},\thinspace V_{\textrm{s}}\right]$
be the orthonormal matrix of eigenvectors of the symmetric primitive
matrix $P_{k}^{T}P_{k}$. By spectral decomposition \cite{Ref:Chung12},
we get: \vspace{-10pt}

\[
V_{\textrm{tr}}^{T}P_{k}^{T}P_{k}V_{\textrm{tr}}=\left[\begin{array}{cc}
1 & \mathbf{0}^{1\times(m-1)}\\
\mathbf{0}^{(m-1)\times1} & V_{\textrm{s}}^{T}P_{k}^{T}P_{k}V_{\textrm{s}}
\end{array}\right]\thinspace,
\]
where $\frac{1}{m}\mathbf{1}^{T}P_{k}^{T}P_{k}\mathbf{1}=1$, $\frac{1}{\sqrt{m}}\mathbf{1}^{T}P_{k}^{T}P_{k}V_{\textrm{s}}=\mathbf{0}^{1\times(m-1)}$,
and $V_{\textrm{s}}^{T}P_{k}^{T}P_{k}\mathbf{1}\frac{1}{\sqrt{m}}=\mathbf{0}^{(m-1)\times1}$
are used. Since the eigenvectors are orthonormal, $V_{s}V_{s}^{T}+\frac{1}{m}\mathbf{1}\mathbf{1}^{T}=\mathbf{I}$.
The rate at which $\mathbf{\mathcal{U}}_{k,\nu}$ synchronizes to
$\frac{1}{\sqrt{m}}\mathbf{1}$ (or $\mathbf{\mathcal{U}}_{k}^{\star}$)
is equal to the rate at which $V_{\textrm{s}}^{T}\mathbf{\mathcal{U}}_{k,\nu}\rightarrow\mathbf{0}^{(m-1)\times1}$.
Pre-multiplying (\ref{log2_PW_eqn2}) by $V_{s}^{T}$ and substituting
$V_{s}^{T}\mathbf{1}=0$ results in: \vspace{-15pt}

\begin{align*}
V_{\textrm{s}}^{T}\mathcal{U}_{k,\nu} & =V_{\textrm{s}}^{T}P_{k}\left(V_{s}V_{s}^{T}+\frac{1}{m}\mathbf{1}\mathbf{1}^{T}\right)\mathcal{U}_{k,\nu-1}\\
 & =V_{\textrm{s}}^{T}P_{k}V_{s}V_{s}^{T}\mathcal{U}_{k,\nu-1}\thinspace.
\end{align*}
Let $\mathbf{z}_{k,\nu}=V_{\textrm{s}}^{T}\mathcal{U}_{k,\nu}$. The
corresponding virtual dynamics is represented by $\mathbf{z}_{k,\nu}=(V_{\textrm{s}}^{T}P_{k}V_{s})\mathbf{z}_{k,\nu-1}$,
which has both $V_{\textrm{s}}^{T}\mathcal{U}_{k,\nu}$ and $\mathbf{0}$
as particular solutions. Let $\Phi_{k,\nu}=\mathbf{z}_{k,\nu}^{T}\mathbf{z}_{k,\nu}$
be a candidate Lyapunov function for this dynamics. Expanding this
gives: \vspace{-15pt}

\begin{align*}
\Phi_{k,\nu} & \!\!=\!\mathbf{z}_{k,\nu-1}^{T}V_{\textrm{s}}^{T}P_{k}^{T}P_{k}V_{s}\mathbf{z}_{k,\nu-1}\!\!\leq\!\left(\!\lambda_{\mathrm{max}}\!(V_{\textrm{s}}^{T}P_{k}^{T}P_{k}V_{s})\!\right)\Phi_{k,\nu-1}.
\end{align*}
Note that $V_{\textrm{s}}^{T}P_{k}^{T}P_{k}V_{s}$ contains all the
eigenvalues of $P_{k}^{T}P_{k}$ other than $1$. Hence $\lambda_{\mathrm{max}}(V_{\textrm{s}}^{T}P_{k}^{T}P_{k}V_{s})=\lambda_{m-1}(P_{k}^{T}P_{k})<1$
and $\Phi_{k,\nu}$ globally exponentially vanishes with a rate faster
or equal to $\lambda_{m-1}(P_{k}^{T}P_{k})$. Hence each $\mathbf{\mathcal{H}}_{k,\nu}^{j}$
globally exponentially converges pointwise to $\mathbf{\mathcal{H}}_{k}^{\star}$
with a rate faster or equal to $\sqrt{\lambda_{m-1}(P_{k}^{T}P_{k})}=\sigma_{m-1}(P_{k})$.

Next, we need to find the rate of convergence of $\mathbf{\mathcal{F}}_{k,\nu}^{j}$
to $\mathbf{\mathcal{F}}_{k}^{\star}$. From the exponential convergence
of $\mathcal{H}_{k,\nu}^{j}$, we get: \vspace{-15pt}

\begin{align}
 & \left|\ln\left[\frac{p_{k,\nu}^{j}(\boldsymbol{x}_{k})}{p_{k}^{\star}(\boldsymbol{x}_{k})}\frac{p_{k}^{\star}(\boldsymbol{\psi}_{k})}{p_{k,\nu}^{j}(\boldsymbol{\psi}_{k})}\right]\right|\nonumber \\
 & \leq\sigma_{m-1}(P_{k})\left|\ln\left[\frac{p_{k,\nu-1}^{j}(\boldsymbol{x}_{k})}{p_{k}^{\star}(\boldsymbol{x}_{k})}\frac{p_{k}^{\star}(\boldsymbol{\psi}_{k})}{p_{k,\nu-1}^{j}(\boldsymbol{\psi}_{k})}\right]\right|\thinspace.\label{simplified_H_to_F_4}
\end{align}
Let us define the function $\alpha_{k,\nu}^{j}(\boldsymbol{x}_{k})$
such that $\alpha_{k,\nu}^{j}(\boldsymbol{x}_{k})=\left[\frac{p_{k,\nu}^{j}(\boldsymbol{x}_{k})}{p_{k}^{\star}(\boldsymbol{x}_{k})}\frac{p_{k}^{\star}(\boldsymbol{\psi}_{k})}{p_{k,\nu}^{j}(\boldsymbol{\psi}_{k})}\right]$
if $p_{k,\nu}^{j}(\boldsymbol{x}_{k})p_{k}^{\star}(\boldsymbol{\psi}_{k})\geq p_{k}^{\star}(\boldsymbol{x}_{k})p_{k,\nu}^{j}(\boldsymbol{\psi}_{k})$
and $\alpha_{k,\nu}^{j}(\boldsymbol{x}_{k})=\left[\frac{p_{k}^{\star}(\boldsymbol{x}_{k})}{p_{k,\nu}^{j}(\boldsymbol{x}_{k})}\frac{p_{k,\nu}^{j}(\boldsymbol{\psi}_{k})}{p_{k}^{\star}(\boldsymbol{\psi}_{k})}\right]$
otherwise. Note that $\alpha_{k,\nu}^{j}(\boldsymbol{x}_{k})$ is
a continuous function since it is a product of continuous functions.
Since $\alpha_{k,\nu}^{j}(\boldsymbol{x}_{k})\geq1$ and $\ln\left(\alpha_{k,\nu}^{j}(\boldsymbol{x}_{k})\right)\geq0,\thinspace\forall\boldsymbol{x}_{k}\in\mathcal{X}$,
(\ref{simplified_H_to_F_4}) simplifies to: \vspace{-15pt}

\begin{align}
\ln\left(\alpha_{k,\nu}^{j}(\boldsymbol{x}_{k})\right) & \leq\sigma_{m-1}(P_{k})\ln\left(\alpha_{k,\nu-1}^{j}(\boldsymbol{x}_{k})\right)\thinspace.\nonumber \\
\alpha_{k,\nu}^{j}(\boldsymbol{x}_{k}) & \leq\left(\alpha_{k,0}^{j}(\boldsymbol{x}_{k})\right)^{\left(\sigma_{m-1}(P_{k})\right)^{\nu}}\thinspace.\label{simplified_H_to_F_6}
\end{align}
Since $p_{k,\nu}^{j}(\boldsymbol{x}_{k})$ tends to $p_{k}^{\star}(\boldsymbol{x}_{k})$,
i.e., $\lim_{\nu\rightarrow\infty}\alpha_{k,\nu}^{j}(\boldsymbol{x}_{k})=1$,
we can write (\ref{simplified_H_to_F_6}) as: \vspace{-10pt}

\begin{equation}
\alpha_{k,\nu}^{j}(\boldsymbol{x}_{k})-1\leq\left(\alpha_{k,0}^{j}(\boldsymbol{x}_{k})\right)^{\left(\sigma_{m-1}(P_{k})\right)^{\nu}}-1^{\left(\sigma_{m-1}(P_{k})\right)^{\nu}}\thinspace.\label{simplified_H_to_F_7}
\end{equation}
Using the mean value theorem (cf. \cite{Ref:Ross80}), the right hand
side of (\ref{simplified_H_to_F_7}) can be simplified to (\ref{simplified_H_to_F_8}),
for some $c\in[1,\alpha_{k,0}^{j}(\boldsymbol{x}_{k})]$. \vspace{-10pt}

\begin{align}
 & \left(\alpha_{k,0}^{j}(\boldsymbol{x}_{k})\right)^{\left(\sigma_{m-1}(P_{k})\right)^{\nu}}-1^{\left(\sigma_{m-1}(P_{k})\right)^{\nu}}\nonumber \\
 & =\left(\sigma_{m-1}(P_{k})\right)^{\nu}\left(c^{\left(\sigma_{m-1}(P_{k})\right)^{\nu}-1}\right)\left(\alpha_{k,0}^{j}(\boldsymbol{x}_{k})-1\right)\thinspace.\label{simplified_H_to_F_8}
\end{align}
As $\sigma_{m-1}(P_{k})<1$, the maximum value of $\left(c^{\left(\sigma_{m-1}(P_{k})\right)^{\nu}-1}\right)$
is $1$. Substituting this result into (\ref{simplified_H_to_F_7})
gives: \vspace{-10pt}

\begin{equation}
\alpha_{k,\nu}^{j}(\boldsymbol{x}_{k})-1\leq\left(\sigma_{m-1}(P_{k})\right)^{\nu}\left(\alpha_{k,0}^{j}(\boldsymbol{x}_{k})-1\right)\thinspace.\label{simplified_H_to_F_9}
\end{equation}
Hence $\alpha_{k,\nu}^{j}(\boldsymbol{x}_{k})$ exponentially converges
to $1$ with a rate faster or equal to $\sigma_{m-1}(P_{k})$. Irrespective
of the orientation of $\alpha_{k,\nu}^{j}(\boldsymbol{x}_{k})$ and
$\alpha_{k,0}^{j}(\boldsymbol{x}_{k})$, (\ref{simplified_H_to_F_9})
can be written as (\ref{simplified_H_to_F_11}) by multiplying with
$\frac{1}{\alpha_{k,\nu}^{j}(\boldsymbol{x}_{k})}$ or $\frac{1}{\alpha_{k,0}^{j}(\boldsymbol{x}_{k})}$,
and then with $p_{k}^{\star}(\boldsymbol{x}_{k})$. \vspace{-10pt}

\begin{align}
 & \left|\frac{p_{k}^{\star}(\boldsymbol{\psi}_{k})}{p_{k,\nu}^{j}(\boldsymbol{\psi}_{k})}p_{k,\nu}^{j}(\boldsymbol{x}_{k})-p_{k}^{\star}(\boldsymbol{x}_{k})\right|\nonumber \\
 & \leq\left(\sigma_{m-1}(P_{k})\right)^{\nu}\left|\frac{p_{k}^{\star}(\boldsymbol{\psi}_{k})}{p_{k,0}^{j}(\boldsymbol{\psi}_{k})}p_{k,0}^{j}(\boldsymbol{x}_{k})-p_{k}^{\star}(\boldsymbol{x}_{k})\right|\thinspace.\label{simplified_H_to_F_11}
\end{align}
As shown in the proof of Theorem \ref{thm:LogOP-SC-digraphs}, we
can choose $\tilde{\boldsymbol{\psi}}_{k}\in\mathcal{X}$ such that
$p_{k,0}^{j}(\tilde{\boldsymbol{\psi}}_{k})=p_{k}^{\star}(\tilde{\boldsymbol{\psi}}_{k})$.
Now we discuss two cases to reduce the left hand side of (\ref{simplified_H_to_F_11})
to $\left|p_{k,\nu}^{j}(\boldsymbol{x}_{k})-p_{k}^{\star}(\boldsymbol{x}_{k})\right|$.
\vspace{-10pt}

\begin{align*}
 & \left|\frac{p_{k}^{\star}(\tilde{\boldsymbol{\psi}}_{k})}{p_{k,\nu}^{j}(\tilde{\boldsymbol{\psi}}_{k})}p_{k,\nu}^{j}(\boldsymbol{x}_{k})-p_{k}^{\star}(\boldsymbol{x}_{k})\right|\\
 & =\begin{cases}
\left|p_{k,\nu}^{j}(\boldsymbol{x}_{k})-p_{k}^{\star}(\boldsymbol{x}_{k})+\left(\frac{p_{k}^{\star}(\tilde{\boldsymbol{\psi}}_{k})}{p_{k,\nu}^{j}(\tilde{\boldsymbol{\psi}}_{k})}-1\right)p_{k,\nu}^{j}(\boldsymbol{x}_{k})\right|\\
\qquad\textrm{ if }\frac{p_{k}^{\star}(\tilde{\boldsymbol{\psi}}_{k})}{p_{k,\nu}^{j}(\tilde{\boldsymbol{\psi}}_{k})}\geq1\\
\left|p_{k}^{\star}(\boldsymbol{x}_{k})-p_{k,\nu}^{j}(\boldsymbol{x}_{k})+\left(1-\frac{p_{k}^{\star}(\tilde{\boldsymbol{\psi}}_{k})}{p_{k,\nu}^{j}(\tilde{\boldsymbol{\psi}}_{k})}\right)p_{k,\nu}^{j}(\boldsymbol{x}_{k})\right|\\
\qquad\textrm{ if }\frac{p_{k}^{\star}(\tilde{\boldsymbol{\psi}}_{k})}{p_{k,\nu}^{j}(\tilde{\boldsymbol{\psi}}_{k})}<1
\end{cases}\\
 & \geq\left|p_{k}^{\star}(\boldsymbol{x}_{k})-p_{k,\nu}^{j}(\boldsymbol{x}_{k})\right|\thinspace.
\end{align*}
Hence, for both the cases, we are able to simplify (\ref{simplified_H_to_F_11})
to: \vspace{-10pt}

\begin{align*}
\left|p_{k,\nu}^{j}(\boldsymbol{x}_{k})-p_{k}^{\star}(\boldsymbol{x}_{k})\right| & \leq\left(\sigma_{m-1}(P_{k})\right)^{\nu}\left|p_{k,0}^{j}(\boldsymbol{x}_{k})-p_{k}^{\star}(\boldsymbol{x}_{k})\right|\thinspace.
\end{align*}
Thus each $\mathcal{F}_{k,\nu}^{j}=p_{k,\nu}^{j}(\boldsymbol{x}_{k})$
globally exponentially converges to $\mathcal{F}_{k}^{\star}=p_{k}^{\star}(\boldsymbol{x}_{k})$
with a rate faster or equal to $\sigma_{m-1}(P_{k})$.
\end{IEEEproof}
The KL divergence is a measure of the information lost when the consensual
pdf is used to approximate the locally estimated posterior pdfs. We
now show that the consensual pdf $\mathcal{F}_{k}^{\star}$ obtained
using Theorem \ref{thm:LogOP-SC-balanced-digraphs}, which is the
weighted geometric average of the locally estimated posterior pdfs
$\mathbf{\mathcal{F}}_{k,0}^{j},\forall j\in\{1,\ldots,m\}$, minimizes
the information lost during the consensus stage because it minimizes
the sum of KL divergences with those pdfs. 

\begin{figure}
\begin{centering}
\begin{tabular}{cc}
\includegraphics[bb=0bp 0bp 364bp 306bp,clip,width=1.65in]{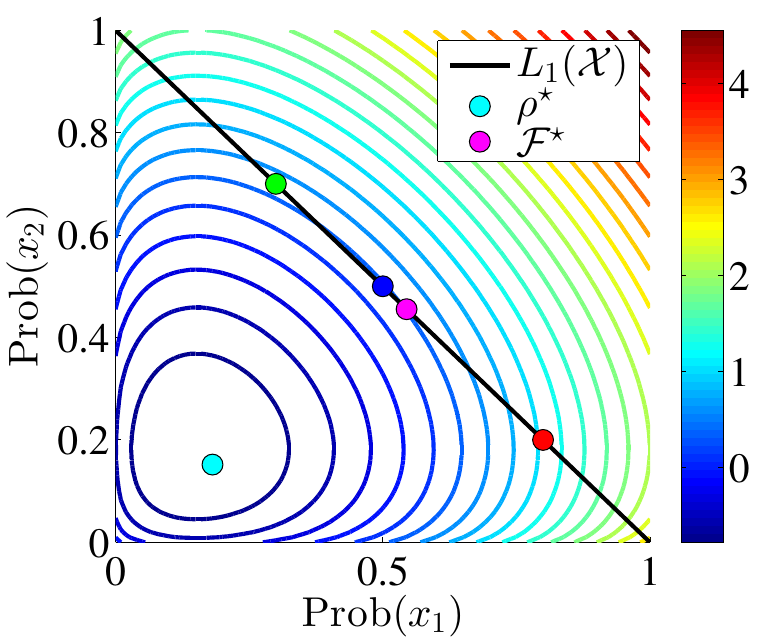} & \includegraphics[bb=0bp 0bp 340bp 306bp,clip,width=1.5in]{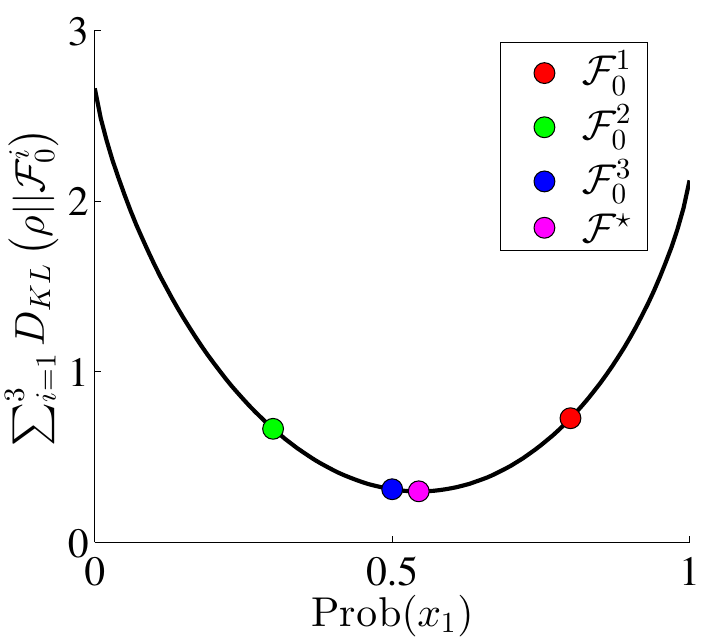}\tabularnewline
(a) & (b)\tabularnewline
\end{tabular}
\par\end{centering}

\caption{Let the discrete state space $\mathcal{X}$ have only two states $x_{1}$
and $x_{2}$. All valid pmfs must lie on the set $\mathscr{L}_{1}(\mathcal{X})$
where $\mathbb{P}(x_{1})+\mathbb{P}(x_{2})=1$. Given three initial
pmfs $\mathcal{F}_{0}^{i},\thinspace i=\{1,2,3\}$, the objective
is to find the pmf that globally minimizes the convex cost function
$\sum_{i=1}^{3}D_{KL}\left(\rho||\mathcal{F}_{0}^{i}\right)$. In
(a), $\rho^{\star}=\frac{1}{e}\prod_{i=1}^{3}(\mathcal{F}_{0}^{i})^{1/m}$
globally minimizes the cost function, but it does not lie on $\mathscr{L}_{1}(\mathcal{X})$.
In (b), $\mathcal{F}^{\star}\in\mathscr{L}_{1}(\mathcal{X})$, which
is the projection of $\rho^{\star}$ on the set $\mathscr{L}_{1}(\mathcal{X})$
obtained by normalizing $\rho^{\star}$ to $1$, indeed globally minimizes
the cost function on the set $\mathscr{L}_{1}(\mathcal{X})$. \label{fig:min-KL-example}}

\end{figure}

\begin{theorem} \label{thm:LogOP-minimizes-KL} The consensual pdf
$\mathcal{F}_{k}^{\star}$ given by (\ref{log2_PW_eqn5}) globally
minimizes the sum of Kullback--Leibler (KL) divergences with the locally
estimated posterior pdfs at the start of the consensus stage $\mathbf{\mathcal{F}}_{k,0}^{j},\forall j\in\{1,\ldots,m\}$,
i.e., \vspace{-10pt}

\begin{equation}
\mathcal{F}_{k}^{\star}=\mathrm{arg}\min_{\rho\in\mathscr{L}_{1}(\mathcal{X})}\thinspace\thinspace\thinspace\sum_{i=1}^{m}D_{KL}\left(\rho||\mathcal{F}_{k,0}^{i}\right),\label{sum_KL_dist_0}
\end{equation}
where $\mathscr{L}_{1}(\mathcal{X})$ is the set of all pdfs over
the state space $\mathcal{X}$ satisfying Assumption \ref{assump:nonnegative_pdf}.
\end{theorem}
\begin{IEEEproof}
The sum of the KL divergences of a pdf $\rho\in\mathscr{L}_{1}(\mathcal{X})$
with the locally estimated posterior pdfs is given by: \vspace{-10pt}

\begin{align}
 & \sum_{i=1}^{m}D_{KL}\left(\rho||\mathcal{F}_{k,0}^{i}\right)=\nonumber \\
 & \sum_{i=1}^{m}\int_{\mathcal{X}}\!\left(\rho(\boldsymbol{x}_{k})\ln(\rho(\boldsymbol{x}_{k}))\!-\!\rho(\boldsymbol{x}_{k})\ln(p_{k,0}^{i}(\boldsymbol{x}_{k}))\right)d\mu(\boldsymbol{x}_{k}).\label{sum_KL_dist_1}
\end{align}
Under Assumption \ref{assump:nonnegative_pdf}, $D_{KL}\left(\rho||\mathcal{F}_{k,0}^{i}\right)$
is well defined for all agents. Differentiating (\ref{sum_KL_dist_1})
with respect to $\rho$ using Leibniz integral rule \cite[Theorem A.5.1, pp. 372]{Ref:Durrett05},
and equating it to zero gives: \vspace{-10pt}

\begin{align}
\sum_{i=1}^{m}\int_{\mathcal{X}}\left(\ln(\rho(\boldsymbol{x}_{k}))+1-\ln(p_{k,0}^{i}(\boldsymbol{x}_{k}))\right)\thinspace d\mu(\boldsymbol{x}_{k})=0\thinspace,\label{sum_KL_dist_2}
\end{align}
where $\rho^{\star}(\boldsymbol{x}_{k})=\frac{1}{e}\prod_{i=1}^{m}(p_{k,0}^{i}(\boldsymbol{x}_{k}))^{1/m}$
is the solution to (\ref{sum_KL_dist_2}). The projection of $\rho^{\star}$
on the set $\mathscr{L}_{1}(\mathcal{X})$, obtained by normalizing
$\rho^{\star}$ to $1$, is the consensual pdf $\mathcal{F}_{k}^{\star}\in\mathscr{L}_{1}(\mathcal{X})$
given by (\ref{log2_PW_eqn5}). 

The KL divergence is a convex function of pdf pairs \cite[Theorem 2.7.2, pp. 30]{Ref:Cover91},
hence the sum of KL divergences (\ref{sum_KL_dist_1}) is a convex
function of $\rho$. If $\rho_{1},\rho_{2},\ldots,\rho_{n}\in\mathscr{L}_{1}(\mathcal{X})$
and $\eta_{1},\eta_{2},\ldots,\eta_{n}\in[0,1]$ such that $\sum_{i=1}^{n}\eta_{i}=1$,
then $\rho^{\dagger}=\sum_{i=1}^{n}\eta_{i}\rho_{i}\in\mathscr{L}_{1}(\mathcal{X})$;
because (i) since $\rho_{i}(\boldsymbol{x}_{k})>0,\forall\boldsymbol{x}_{k}\in\mathcal{X},\forall i\in\{1,\ldots,n\}$
therefore $\rho^{\dagger}(\boldsymbol{x}_{k})>0,\forall\boldsymbol{x}_{k}\in\mathcal{X}$;
and (ii) since $\int_{\mathcal{X}}\rho_{i}(\boldsymbol{x}_{k})d\mu(\boldsymbol{x}_{k})=1,\forall i\in\{1,\ldots,n\}$
therefore $\int_{\mathcal{X}}\rho^{\dagger}(\boldsymbol{x}_{k})d\mu(\boldsymbol{x}_{k})=1$.
Moreover, since $\mathcal{X}$ is a compact set, therefore $\mathscr{L}_{1}(\mathcal{X})$
is a closed set. Hence $\mathscr{L}_{1}(\mathcal{X})$ is a closed
convex set. Hence (\ref{sum_KL_dist_0}) is a convex optimization
problem. 

The gradient of $\sum_{i=1}^{m}D_{KL}\left(\rho||\mathcal{F}_{k,0}^{i}\right)$
evaluated at $\mathcal{F}_{k}^{\star}$ is a constant, i.e., \vspace{-15pt}

\begin{align*}
 & \left.\frac{d}{d\rho}\!\sum_{i=1}^{m}\! D_{KL}\!\left(\rho||\mathcal{F}_{k,0}^{i}\right)\!\right|_{\textrm{\!}\rho=\mathcal{F}_{k}^{\star}}\!\!\!\!\!\!=\! m\ln\frac{e}{\int_{\mathcal{X}}\!\prod_{i=1}^{m}\!\left(p_{k,0}^{i}(\boldsymbol{x}_{k})\right)^{\frac{1}{m}}\! d\mu(\boldsymbol{x}_{k})}.
\end{align*}
This indicates that for further minimizing the convex cost function,
we have to change the normalizing constant of $\mathcal{F}_{k}^{\star}$,
which will result in exiting the set $\mathscr{L}_{1}(\mathcal{X})$.
Hence $\mathcal{F}_{k}^{\star}$ is the global minimum of the convex
cost function (\ref{sum_KL_dist_0}) in the convex set $\mathscr{L}_{1}(\mathcal{X})$.
This is illustrated using a simple example in Fig. \ref{fig:min-KL-example}. 

Another proof approach involves taking the logarithm, in the KL divergence
formula, to the base $c:=\left(\int_{\mathcal{X}}\!\prod_{i=1}^{m}\!\left(p_{k,0}^{i}(\boldsymbol{x}_{k})\right)^{\frac{1}{m}}\! d\mu(\boldsymbol{x}_{k})\right)$.
Then differentiating $\sum_{i=1}^{m}D_{KL}\left(\rho||\mathcal{F}_{k,0}^{i}\right)$
with respect to $\rho$ gives: \vspace{-10pt}

\[
\sum_{i=1}^{m}\int_{\mathcal{X}}\left(\log_{c}(\rho(\boldsymbol{x}_{k}))+1-\log_{c}(p_{k,0}^{i}(\boldsymbol{x}_{k}))\right)\thinspace d\mu(\boldsymbol{x}_{k})=0\thinspace,
\]
which is minimized by $\mathcal{F}_{k}^{\star}$. Hence $\mathcal{F}_{k}^{\star}$
is indeed the global minimum of the convex optimization problem (\ref{sum_KL_dist_0}). 
\end{IEEEproof}
Note that if a central agent receives all the locally estimated posterior
pdfs ($\mathbf{\mathcal{F}}_{k,0}^{j},\forall j\in\{1,\ldots,m\}$)
and is tasked to find the best estimate in the information theoretic
sense, then it would also yield the same consensual pdf $\mathcal{F}_{k}^{\star}$
given by (\ref{log2_PW_eqn5}). Hence we claim to have achieved distributed
estimation using this algorithm. 

In Remark \ref{remark:gQLOP}, we state that the methods for recursively
combining probability distributions to reach a consensual distribution
are limited to LinOP, LogOP, and their affine combinations. 

\begin{remark} \label{remark:gQLOP} The LinOP and LogOP methods
for combining probability distributions can be generalized by the
g--Quasi--Linear Opinion Pool (g--QLOP), which is described by the
following equation: \vspace{-15pt}

\begin{align}
\mathbf{\mathcal{F}}_{k,\nu}^{j}= & \frac{g^{-1}(\sum_{\ell\in\mathcal{J}_{k}^{j}}\alpha_{k,\nu-1}^{j\ell}g(\mathbf{\mathcal{F}}_{k,\nu-1}^{\ell}))}{\int_{\mathcal{X}}g^{-1}(\sum_{\ell\in\mathcal{J}_{k}^{j}}\alpha_{k,\nu-1}^{j\ell}g(\mathbf{\mathcal{F}}_{k,\nu-1}^{\ell}))\thinspace d\mu(\boldsymbol{x}_{k})},\nonumber \\
 & \forall j\in\{1,\ldots,m\},\forall\nu\in\mathbb{N},\label{gqlop_agreement_equation}
\end{align}
where $g$ is a continuous, strictly monotone function. It is shown
in \cite{Ref:Gilardoni93} that, other than the linear combination
of LinOP and LogOP, there is no function $g$ for which the final
consensus can be expressed by the following equation: \vspace{-10pt}

\begin{align}
\lim_{\nu\rightarrow\infty}\mathcal{F}_{k,\nu}^{j}= & \mathbf{\mathcal{F}}_{k}^{\star}=\frac{g^{-1}(\sum_{j=1}^{m}\pi_{j}g(\mathcal{F}_{k,0}^{\ell}))}{\int_{\mathcal{X}}g^{-1}(\sum_{j=1}^{m}\pi_{j}g(\mathcal{F}_{k,0}^{\ell}))\thinspace d\mu(\boldsymbol{x}_{k})},\nonumber \\
 & \forall j\in\{1,\ldots,m\},\forall\nu\in\mathbb{N},\label{gqlop_consensus}
\end{align}
where $\boldsymbol{\pi}$ is the unique stationary solution. Moreover,
the function $g$ is said to be k-Markovian if the scheme for combining
probability distribution (\ref{gqlop_agreement_equation}) yields
the consensus (\ref{gqlop_consensus}) for every regular communication
network topology and for all initial positive densities. It is also
shown that $g$ is k-Markovian if and only if the g--QLOP is either
 LinOP or LogOP \cite{Ref:Gilardoni93}. \hfill $\Box$  \end{remark}

\subsection{Communicating Probability Distributions \label{sub:Communicating-probability-distri}}

The consensus algorithms using either LinOP or LogOP need the estimated
pdfs to be communicated to other agents in the network. We propose
to adopt the following methods for communicating pdfs.

The first approach involves approximating the pdf by a weighted sum
of Gaussians and then transmitting this approximate distribution.
Let $\mathcal{N}(\boldsymbol{x}_{k}-m_{i},B_{i})$ denote the Gaussian
density function, where the mean is the $n_{x}$-vector $m_{i}$ and
the covariance is the positive-definite symmetric matrix $B_{i}$.
The Gaussian sum approximations lemma of \cite[pp. 213]{Ref:Anderson05}
states that any pdf $\mathbf{\mathcal{F}}=p(\boldsymbol{x}_{k})$
can be approximated as closely as desired in the $\mathscr{L}_{1}(\mathbb{R}^{n_{x}})$
space by a pdf of the form $\mathbf{\hat{\mathcal{F}}}=\hat{p}(\boldsymbol{x}_{k})=\sum_{i=1}^{n_{g}}\alpha_{i}\mathcal{N}(\boldsymbol{x}_{k}-m_{i},B_{i})$,
for some integer $n_{g}$ and positive scalars $\alpha_{i}$ with
$\sum_{i=1}^{n_{g}}\alpha_{i}=1$. For an acceptable communication
error $\varepsilon_{\textrm{comm}}>0$, there exists $n_{g}$, $\alpha_{i}$,
$m_{i}$ and $B_{i}$ such that $\|\mathbf{\mathcal{F}}-\hat{\mathbf{\mathcal{F}}}\|_{\mathcal{L}_{1}}\leq\varepsilon_{\textrm{comm}}$.
Several techniques for estimating the parameters are discussed in
the Gaussian mixture model literature, like maximum likelihood (ML)
and maximum a posteriori (MAP) parameter estimation \nocite{Ref:Reynolds08,Ref:McLachlan88,Ref:Kotecha03}
\cite{Ref:Reynolds08}--\cite{Ref:Kotecha03}. Hence, in order to
communicate the pdf $\mathcal{\hat{F}}$, the agent needs to transmit
$\frac{1}{2}n_{g}n_{x}\left(n_{x}+3\right)$ real numbers. 

Let us study the effect of this communication error $\varepsilon_{\textrm{comm}}$
on the LinOP consensual distribution. Let $\tilde{\mathcal{F}}_{k,\nu}^{j}$
be the LinOP solution after combining local pdfs corrupted by communication
error, i.e., $\tilde{\mathcal{F}}_{k,\nu}^{j}:=\mathcal{T}\left(\cup_{\ell\in\mathcal{J}_{k}^{j}}\{\hat{\mathcal{F}}_{k,\nu-1}^{\ell}\}\right)$
where $\mathcal{T}(\cdot)$ is LinOP (\ref{agreement_equation}).
We prove by induction that $\|\mathbf{\mathcal{F}}_{k,\nu}^{j}-\tilde{\mathcal{F}}_{k,\nu}^{j}\|_{\mathcal{L}_{1}}\leq\nu\varepsilon_{\mathrm{comm}},\thinspace\forall\nu\in\mathbb{N}$,
where $\mathbf{\mathcal{F}}_{k,\nu}^{j}$ is the true solution obtained
from uncorrupted local pdfs. As the basis of induction holds, the
inductive step for the $\nu^{\textrm{th}}$ consensus step is as follows:
\vspace{-10pt}

\begin{align}
 & \|\mathbf{\mathcal{F}}_{k,\nu}^{j}-\tilde{\mathcal{F}}_{k,\nu}^{j}\|_{\mathcal{L}_{1}}\leq\sum_{\ell\in\mathcal{J}_{k}^{j}}a_{k,\nu-1}^{j\ell}\left(\|\mathbf{\mathcal{F}}_{k,\nu-1}^{\ell}-\tilde{\mathcal{F}}_{k,\nu-1}^{\ell}\|_{\mathcal{L}_{1}}\right.\nonumber \\
 & \quad\left.+\|\tilde{\mathcal{F}}_{k,\nu-1}^{\ell}-\hat{\mathbf{\mathcal{F}}}_{k,\nu-1}^{\ell}\|_{\mathcal{L}_{1}}\right)\leq(\nu-1)\varepsilon_{\textrm{comm}}+\varepsilon_{\textrm{comm}}.\label{eq:comm-error}
\end{align}
Similarly, it follows from the proof of Theorem \ref{thm:LogOP-SC-balanced-digraphs}
that the LogOP solution after $n_{\textrm{loop}}$ iterations ($\tilde{\mathcal{F}}_{k,n_{\textrm{loop}}}^{j}$),
under communication inaccuracies, is always within a ball of $n_{\textrm{loop}}\varepsilon_{\textrm{comm}}$
radius from the true solution using LogOP ($\mathcal{F}_{k,n_{\textrm{loop}}}^{j}$)
in the $\mathscr{L}_{1}(\mathcal{X})$ space, i.e., $\|\mathcal{F}_{k,n_{\textrm{loop}}}^{j}-\tilde{\mathcal{F}}_{k,n_{\textrm{loop}}}^{j}\|_{\mathcal{L}_{1}}\leq n_{\textrm{loop}}\varepsilon_{\textrm{comm}}$. 

If particle filters are used to evaluate the Bayesian filter and combine
the pdfs \cite{Ref:Arulampalam02,Ref:Kotecha03}, then the resampled
particles represent the agent's estimated pdf of the target. Hence
communicating pdfs is equivalent to transmitting these resampled particles. 

The information theoretic approach for communicating pdfs is discussed
in \cite{Ref:Kramer07}. Let the local pdf $\mathcal{F}_{k,\nu}^{j}$
be transmitted over a communication channel using a finite sequence
and the pdf $\mathcal{\hat{F}}_{k,\nu}^{j}$ is reconstructed by the
other agent. For a given error threshold, the minimum rate such that
the variational distortion between $\mathcal{F}_{k,\nu}^{j}$ and
$\mathcal{\hat{F}}_{k,\nu}^{j}$ is bounded by the error threshold,
is given by the mutual information between transmitted and received
finite sequences.

Now that we have established that communication of pdfs is possible,
let us discuss the complete BCF algorithm. 

\begin{algorithm}[th]
\caption{BCF--LogOP on SC Balanced Digraphs \label{alg:Dynamic-Bayesian-Consensus}}

\begin{singlespace}
{\small{}}%
\begin{tabular}{rll}
{\small{1:}} & \multicolumn{2}{l}{{\small{(one cycle of $j^{\textrm{th}}$ agent during $k^{\textrm{th}}$
time instant)}}}\tabularnewline
{\small{2:}} & \multicolumn{2}{l}{{\small{Given the pdf from previous time step }}}\tabularnewline
 & {\small{$\:\:\:$ $\mathbf{\mathcal{F}}_{k-1}^{j}=p_{k-1}^{j}(\boldsymbol{x}_{k-1})$}} & \tabularnewline
{\small{3:}} & {\small{Set $n_{\textrm{loop}}$, the weights $a_{k}^{j\ell}$ }} & $\}$ Theorems \ref{thm:LogOP-SC-balanced-digraphs}, \ref{thm:BCF-LogOP-nloop}\tabularnewline
{\small{4:}} & \textbf{\small{while }}{\small{tracking }}\textbf{\small{do}} & \tabularnewline
{\small{5:}} & {\small{$\:\:\:$Compute the prior pdf }} & \multirow{5}{*}{$\left\} \begin{array}{c}
\textrm{Bayesian}\\
\textrm{Filtering}\\
\textrm{Stage}\\
\textrm{(Sec. \ref{sub:Bayesian-Filtering-Algorithms})}\\
\\
\end{array}\right.$}\tabularnewline
 & \textbf{\small{$\:\:\:\:\:\:$}}{\small{$p_{k}^{j}(\boldsymbol{x}_{k})$
using (\ref{predict_stage_con})}} & \tabularnewline
{\small{6:}} & \textbf{\small{$\:\:\:$}}{\small{Compute the posterior pdf }} & \tabularnewline
 & \textbf{\small{$\:\:\:$$\:\:\:$}}{\small{$\mathbf{\mathcal{F}}_{k}^{j}=p_{k}^{j}(\boldsymbol{x}_{k}|\boldsymbol{z}_{k}^{\mathcal{S}_{k}^{j}})$
using (\ref{update_stage1_con})}} & \tabularnewline
 & \textbf{\small{$\:\:\:$$\:\:\:$ }}{\small{and measurement array
$\boldsymbol{z}_{k}^{\mathcal{S}_{k}^{j}}$ }} & \tabularnewline
{\small{7:}} & \textbf{\small{$\:\:\:$for }}{\small{$\nu=1$ to $n_{\textrm{loop}}$}} & \multirow{8}{*}{$\left\} \begin{array}{c}
\\
\\
\\
\textrm{LogOP--based}\\
\textrm{Consensus}\\
\textrm{Stage}\\
\textrm{(Sec. \ref{sub:Logarithmic-Opinion-Pool})}\\
\\
\\
\\
\end{array}\right.$}\tabularnewline
{\small{8:}} & \textbf{\small{$\:\:\:\:\:\:$if $\nu=1$ then }}{\small{Set $\mathbf{\mathcal{F}}_{k,0}^{j}=\mathbf{\mathcal{F}}_{k}^{j}$}} & \tabularnewline
 & \textbf{\small{$\:\:\:\:\:\:$end if}} & \tabularnewline
{\small{9:}} & \textbf{\small{$\:\:\:\:\:\:$}}{\small{Obtain the communicated }} & \tabularnewline
 & \textbf{\small{$\:\:\:\:\:\:$$\:\:\:$}}{\small{pdfs $\mathbf{\mathcal{F}}_{k,\nu-1}^{\ell},\forall\ell\in\mathcal{J}_{k}^{j}$}} & \tabularnewline
{\small{10:}} & \textbf{\small{$\:\:\:\:\:\:$}}{\small{Compute the new pdf $\mathbf{\mathcal{F}}_{k,\nu}^{j}$ }} & \tabularnewline
 & \textbf{\small{$\:\:\:\:\:\:$$\:\:\:$}}{\small{using the LogOP (\ref{log_agreement_equation})}} & \tabularnewline
 & \textbf{\small{$\:\:\:$end for}} & \tabularnewline
{\small{11:}} & \textbf{\small{$\:\:\:$}}{\small{Set $\mathbf{\mathcal{F}}_{k}^{j}=\mathbf{\mathcal{F}}_{k,n_{\textrm{loop}}}^{j}$}} & \tabularnewline
 & \textbf{\small{end while}} & \tabularnewline
\end{tabular}\end{singlespace}
\end{algorithm}

\vspace{-10pt}

\section{Main Algorithms: Bayesian Consensus Filtering \label{sec:Dynamic-Bayesian-Consensus}}

In this section, we finally solve the complete problem statement for
BCF discussed in Section \ref{sub:Problem-Statement} and \textbf{Algorithm
\ref{alg:Dynamic-Bayesian-Consensus}}. We also introduce an hierarchical
algorithm that can be used when some agents in the network fail to
observe the target. \vspace{-10pt}

\subsection{Bayesian Consensus Filtering \label{sub:Bayesian-Consensus-Filtering}}

The BCF is performed in two steps: (i) each agent locally estimates
the pdf of the target's states using a Bayesian filter with/without
measurements from neighboring agents, as discussed in Section \ref{sub:Bayesian-Filtering-Algorithms},
and (ii) during the consensus stage, each agent recursively transmits
its pdf estimate of the target's states to other agents, receives
estimates of its neighboring agents, and combines them using the LogOP
as discussed in Section \ref{sub:Logarithmic-Opinion-Pool}. According
to \cite{Ref:Sayed12}, this strategy of first updating the local
estimate and then combining these local estimates to achieve a consensus
is stable and gives the best performance in comparison with other
update--combine strategies. In this section, we compute the number
of consensus loops ($n_{\textrm{loop}}$ in \textbf{Algorithm \ref{alg:Dynamic-Bayesian-Consensus}})
needed to reach a satisfactory consensus estimate across the network
and discuss the convergence of this algorithm. 

\begin{definition} (\textit{Disagreement vector $\boldsymbol{\theta}_{k,\nu}$})
Let us define $\boldsymbol{\theta}_{k,\nu}:=\left(\theta_{k,\nu}^{1},\ldots,\theta_{k,\nu}^{m}\right)^{T}$,
where $\theta_{k,\nu}^{j}:=\|\mathbf{\mathcal{F}}_{k,\nu}^{j}-\mathcal{F}_{k}^{\star}\|_{\mathcal{L}_{1}}$.
Since the $\mathcal{L}_{1}$ distance between pdfs is upper bounded
by $2$, the $\ell_{2}$ norm of the disagreement vector ($\|\boldsymbol{\theta}_{k,\nu}\|_{\ell_{2}}$)
is upper bounded by $2\sqrt{m}$. \hfill $\Box$ \end{definition}

This conservative bound is used to obtain the minimum number of consensus
loops for achieving $\varepsilon$-consensus across the network, while
tracking a moving target. Let us now quantify the divergence of the
local pdfs during the Bayesian filtering stage. 

\begin{definition} (\textit{Error propagation dynamics $\boldsymbol{\Gamma}(\cdot)$})
Let us assume that the dynamics of the $\ell_{2}$ norm of the disagreement
vector during the Bayesian filtering stage can be obtained from the
target dynamics and measurement models (\ref{sys_mod}) and (\ref{mes_mod_con}).
The error propagation dynamics $\boldsymbol{\Gamma}(\cdot)$ estimates
the maximum divergence of the local pdfs during the Bayesian filtering
stage, i.e., $\|\boldsymbol{\theta}_{k,0}\|_{\ell_{2}}\leq\boldsymbol{\Gamma}\left(\|\boldsymbol{\theta}_{k-1,n_{\textrm{loop}}}\|_{\ell_{2}}\right)$,
where $\|\boldsymbol{\theta}_{k-1,n_{\textrm{loop}}}\|_{\ell_{2}}$
is the disagreement vector with respect to $\mathcal{F}_{k-1}^{\star}$
at the end of the consensus stage during the $(k-1)^{\textrm{th}}$
time instant; and $\|\boldsymbol{\theta}_{k,0}\|_{\ell_{2}}$ is the
disagreement vector with respect to $\mathcal{F}_{k}^{\star}$ after
the update stage during the $k^{\textrm{th}}$ time instant. \hfill
$\Box$  \end{definition} 

Next we obtain the minimum number of consensus loops for achieving
$\varepsilon$-consensus across the network and also derive conditions
on the communication network topology for a given number of consensus
loops. 

\begin{theorem} \label{thm:BCF-LogOP-nloop}\emph{ (BCF--LogOP on
SC Balanced Digraphs)} Under Assumptions \ref{assump:prior}, \ref{assump:nonnegative_pdf},
\ref{assump:weights2}, and an acceptable communication error $\varepsilon_{\mathrm{comm}}>0$,
each agent tracks the target using the BCF algorithm. For some acceptable
consensus error $\varepsilon_{\mathrm{cons}}>0$ and $\gamma_{k}=\min\left(\boldsymbol{\Gamma}\left(\|\boldsymbol{\theta}_{k-1,n_{\mathrm{loop}}}\|_{\ell_{2}}\right),2\sqrt{m}\right)$:
\\* (i) for a given $P_{k}$, if the number of consensus loops $n_{\mathrm{loop}}$
satisfies \vspace{-5pt}

\begin{equation}
\left(\sigma_{m-1}(P_{k})\right)^{n_{\mathrm{loop}}}\gamma_{k}+2n_{\mathrm{loop}}\varepsilon_{\mathrm{comm}}\sqrt{m}\leq\varepsilon_{\mathrm{cons}}\thinspace;\label{eq:nloop-givn-sigma-Pk}
\end{equation}
or (ii) for a given $n_{\mathrm{loop}}$, if the communication network
topology ($P_{k}$) during the $k^{\textrm{th}}$ time instant is
such that \vspace{-5pt}

\begin{align}
\sigma_{m-1}(P_{k})\leq\left(\frac{\varepsilon_{\mathrm{cons}}-2n_{\mathrm{loop}}\varepsilon_{\mathrm{comm}}\sqrt{m}}{\gamma_{k}}\right)^{\frac{1}{n_{\mathrm{loop}}}}\thinspace;\label{eq:sigma-Pk-given-nloop}
\end{align}
then the $\ell_{2}$ norm of the disagreement vector at the end of
the consensus stage is less than $\varepsilon_{\mathrm{cons}}$, i.e.,
$\|\boldsymbol{\theta}_{k,n_{\mathrm{loop}}}\|_{\ell_{2}}\leq\varepsilon_{\mathrm{cons}}$.
\end{theorem} 
\begin{IEEEproof}
In the absence of communication inaccuracies, Theorem \ref{thm:LogOP-SC-balanced-digraphs}
states that the local estimated pdfs $\mathcal{F}_{k,0}^{j}$ globally
exponentially converges pointwise to a consensual pdf $\mathcal{F}_{k}^{\star}$
given by (\ref{log2_PW_eqn5}) with a rate of $\sigma_{m-1}(P_{k})$,
i.e. $\|\mathbf{\mathcal{F}}_{k,\nu}^{j}-\mathcal{F}_{k}^{\star}\|_{\mathcal{L}_{1}}\leq\left(\sigma_{m-1}(P_{k})\right)^{\nu}\|\mathbf{\mathcal{F}}_{k,0}^{j}-\mathcal{F}_{k}^{\star}\|_{\mathcal{L}_{1}}$.
If $\boldsymbol{\theta}_{k,0}$ is the initial disagreement vector
at the start of the consensus stage, then $\|\boldsymbol{\theta}_{k,n_{\textrm{loop}}}\|_{\ell_{2}}\leq\left(\sigma_{m-1}(P_{k})\right)^{n_{\textrm{loop}}}\|\boldsymbol{\theta}_{k,0}\|_{\ell_{2}}\leq\left(\sigma_{m-1}(P_{k})\right)^{n_{\textrm{loop}}}\gamma_{k}$. 

In the presence of communication error, combining (\ref{eq:comm-error})
with the previous result gives $\|\tilde{\mathbf{\mathcal{F}}}_{k,\nu}^{j}-\mathcal{F}_{k}^{\star}\|_{\mathcal{L}_{1}}\leq\left(\sigma_{m-1}(P_{k})\right)^{\nu}\|\mathbf{\mathcal{F}}_{k,0}^{j}-\mathcal{F}_{k}^{\star}\|_{\mathcal{L}_{1}}+\nu\varepsilon_{\textrm{comm}}$.
Since $\theta_{k,\nu}^{j}\leq\|\tilde{\mathbf{\mathcal{F}}}_{k,\nu}^{j}-\mathcal{F}_{k}^{\star}\|_{\mathcal{L}_{1}}+\nu\varepsilon_{\textrm{comm}}$,
the disagreement vector after $n_{\textrm{loop}}$ iterations is given
by $\|\boldsymbol{\theta}_{k,n_{\textrm{loop}}}\|_{\ell_{2}}\leq\left(\sigma_{m-1}(P_{k})\right)^{n_{\textrm{loop}}}\|\boldsymbol{\theta}_{k,0}\|_{\ell_{2}}+2n_{\textrm{loop}}\varepsilon_{\textrm{comm}}\sqrt{m}$.
Thus, we get the conditions on $n_{\textrm{loop}}$ or $\sigma_{m-1}(P_{k})$
from the inequality $\left(\sigma_{m-1}(P_{k})\right)^{n_{\textrm{loop}}}\gamma_{k}+2n_{\textrm{loop}}\varepsilon_{\textrm{comm}}\sqrt{m}\leq\varepsilon_{\mathrm{cons}}$. 
\end{IEEEproof}
Note that in the absence of communication inaccuracies, (\ref{eq:nloop-givn-sigma-Pk})
simplifies to $n_{\mathrm{loop}}\geq\left\lceil \frac{\ln\left(\varepsilon_{\mathrm{cons}}/\gamma_{k}\right)}{\ln\sigma_{m-1}(P_{k})}\right\rceil $
and (\ref{eq:sigma-Pk-given-nloop}) simplifies to $\sigma_{m-1}(P_{k})\leq\left(\frac{\varepsilon_{\mathrm{cons}}}{\gamma_{k}}\right)^{\frac{1}{n_{\mathrm{loop}}}}$.
In the particular case where $n_{\textrm{loop}}=1$ and communication
errors are present, (\ref{eq:sigma-Pk-given-nloop}) simplifies to
$\sigma_{m-1}(P_{k})\leq\frac{\varepsilon_{\mathrm{cons}}-2\varepsilon_{\textrm{comm}}\sqrt{m}}{\gamma_{k}}$
and the necessary condition for a valid solution is $2\varepsilon_{\textrm{comm}}\sqrt{m}<\varepsilon_{\mathrm{cons}}$.
In the genral case, it is desireable that $2\varepsilon_{\mathrm{comm}}\sqrt{m}\ll\varepsilon_{\mathrm{cons}}$
for a valid solution to Theorem \ref{thm:BCF-LogOP-nloop}.

\subsection{Hierarchical Bayesian Consensus Filtering \label{sub:Hierarchical-Bayesian-Consensus}}

In this section, we modify the original problem statement such that
only $m_{1}$ out of $m$ agents are able to observe the target at
the $k^{\textrm{th}}$ time instant. In this scenario, the other $m_{2}(=m-m_{1})$
agents are not able to observe the target. Without loss of generality,
we assume that the first $m_{1}$ agents, i.e., $\{1,2,\ldots,m_{1}\}$,
are tracking the target. During the Bayesian filtering stage, each
\textit{tracking agent} (i.e., agent tracking the target) estimates
the posterior pdf of the target's states at the $k^{\textrm{th}}$
time instant ($\mathbf{\mathcal{F}}_{k}^{j}=p_{k}^{j}(\boldsymbol{x}_{k}|\boldsymbol{z}_{k}^{\mathcal{S}_{k}^{j}\cap\{1,\ldots m_{1}\}}),\forall j\in\{1,\ldots,m_{1}\}$)
using the estimated prior pdf of the target's states ($\mathbf{\mathcal{F}}_{k-1}^{j}$)
and the new measurement array $\boldsymbol{z}_{k}^{\mathcal{S}_{k}^{j}\cap\{1,\ldots m_{1}\}}:=\left\{ \boldsymbol{z}_{k}^{\ell},\forall\ell\in\mathcal{S}_{k}^{j}\cap\{1,\ldots m_{1}\}\right\} $
obtained from the neighboring tracking agents. Each \textit{non-tracking
agent} (i.e., agent not tracking the target) only propagates its prior
pdf during this stage to obtain $p_{k}^{j}(\boldsymbol{x}_{k}),\forall j\in\{m_{1}+1,\ldots,m\}$. 

The objective of hierarchical consensus algorithm is to guarantee
pointwise convergence of each $\mathbf{\mathcal{F}}_{k,\nu}^{j},\forall j\in\{1,\ldots,m\}$
to a pdf $\mathbf{\mathcal{F}}_{k}^{\star}$and only the local estimates
of the agents tracking the target contribute to the consensual pdf.
This is achieved by each tracking agent recursively transmitting its
estimate of the target's states to other agents, only receiving estimates
from its neighboring tracking agents and updating its estimate of
the target. On the other hand, each non-tracking agent recursively
transmits its estimate of the target's states to other agents, receives
estimates from all its neighboring agents and updates its estimate
of the target. This is illustrated using the pseudo-code in \textbf{Algorithm
\ref{alg:HBCF-LogOP}} and the following equations:  \vspace{-15pt}

\begin{align}
\mathbf{\mathcal{F}}_{k,\nu}^{j} & =\mathcal{T}\left(\cup_{\ell\in\mathcal{J}_{k}^{j}\cap\{1,\ldots,m_{1}\}}\{\mathbf{\mathcal{F}}_{k,\nu-1}^{\ell}\}\right),\nonumber \\
 & \forall j\in\{1,\ldots,m_{1}\},\thinspace\forall\nu\in\mathbb{N},\label{combining_estimates_hierachical1}\\
\mathbf{\mathcal{F}}_{k,\nu}^{j} & =\mathcal{T}\left(\cup_{\ell\in\mathcal{J}_{k}^{j}}\{\mathbf{\mathcal{F}}_{k,\nu-1}^{\ell}\}\right),\nonumber \\
 & \forall j\in\{m_{1}+1,\ldots,m\},\thinspace\forall\nu\in\mathbb{N},\label{combining_estimates_hierachical2}
\end{align}
where, $\mathcal{T}(\cdot)$ refers to the LogOP (\ref{log_agreement_equation})
for combining pdfs. Let $\mathcal{D}_{k}$ represent the communication
network topology of only the tracking agents. 

\begin{algorithm}[th]
\caption{Hierarchical BCF--LogOP on SC Balanced Digraphs \label{alg:HBCF-LogOP}}

\begin{singlespace}
{\small{}}%
\begin{tabular}{rll}
{\small{1:}} & \multicolumn{2}{l}{{\small{(one cycle of $j^{\textrm{th}}$ agent during $k^{\textrm{th}}$
time instant)}}}\tabularnewline
{\small{2:}} & \multicolumn{2}{l}{{\small{Given the pdf from previous time step }}}\tabularnewline
 & {\small{$\:\:\:$ $\mathbf{\mathcal{F}}_{k-1}^{j}=p_{k-1}^{j}(\boldsymbol{x}_{k-1})$}} & \tabularnewline
{\small{3:}} & {\small{Set $n_{\textrm{loop}}$, the weights $a_{k}^{j\ell}$ }} & $\}$ Theorems \ref{thm:BCF-LogOP-nloop}, \ref{thm:Hierarchical-LogOP-SC-balanced-digraphs}\tabularnewline
{\small{4:}} & \textbf{\small{while }}{\small{tracking }}\textbf{\small{do}} & \tabularnewline
{\small{5:}} & {\small{$\:\:\:$ Compute prior pdf $p_{k}^{j}(\boldsymbol{x}_{k})$
using (\ref{predict_stage_con})}} & \multirow{5}{*}{$\left\} \begin{array}{c}
\textrm{Bayesian}\\
\textrm{Filtering}\\
\textrm{Stage}\\
\textrm{(Sec. \ref{sub:Bayesian-Filtering-Algorithms})}\\
\\
\end{array}\right.$}\tabularnewline
{\small{6:}} & {\small{$\:\:\:$}}\textbf{\small{ if $j\leq m_{1}$ then}} & \tabularnewline
{\small{7:}} & \textbf{\small{$\:\:\:$}}{\small{$\:\:\:$ Compute the posterior
pdf $\mathbf{\mathcal{F}}_{k}^{j}$ }} & \tabularnewline
 & \textbf{\small{$\:\:\:$$\:\:\:$$\:\:\:$}}{\small{ using (\ref{update_stage1_con})
and $\boldsymbol{z}_{k}^{\mathcal{S}_{k}^{j}\cap\{1,\ldots,m_{1}\}}$ }} & \tabularnewline
 & {\small{$\:\:\:$}}\textbf{\small{ end if}} & \tabularnewline
{\small{8:}} & \textbf{\small{$\:\:\:$ for }}{\small{$\nu=1$ to $n_{\textrm{loop}}$}} & \multirow{15}{*}{$\left\} \begin{array}{c}
\\
\\
\\
\\
\\
\textrm{Hierarchical}\\
\textrm{LogOP--based}\\
\textrm{Consensus}\\
\textrm{Stage}\\
\textrm{(Sec. \ref{sub:Hierarchical-Bayesian-Consensus})}\\
\\
\\
\\
\\
\\
\end{array}\right.$}\tabularnewline
{\small{9:}} & \textbf{\small{$\:\:\:\:\:\:$ if $\nu=1$ then }} & \tabularnewline
{\small{10:}} & \textbf{\small{$\:\:\:\:\:\:$$\:\:\:$ if $j\leq m_{1}$ then }}{\small{Set
$\mathbf{\mathcal{F}}_{k,0}^{j}=\mathbf{\mathcal{F}}_{k}^{j}$}} & \tabularnewline
{\small{11:}} & \textbf{\small{$\:\:\:\:\:\:$$\:\:\:$ else}}{\small{ Set $\mathbf{\mathcal{F}}_{k,0}^{j}=p_{k}^{j}(\boldsymbol{x}_{k})$
}}\textbf{\small{end if}} & \tabularnewline
 & \textbf{\small{$\:\:\:\:\:\:$ end if}} & \tabularnewline
{\small{12:}} & \textbf{\small{$\:\:\:\:\:\:$ if $j\leq m_{1}$ then}} & \tabularnewline
{\small{13:}} & \textbf{\small{$\:\:\:\:\:\:\:\:\:$}}{\small{ Obtain the pdfs $\mathbf{\mathcal{F}}_{k,\nu-1}^{\ell},$ }} & \tabularnewline
 & \textbf{\small{$\:\:\:\:\:\:\:\:\:$$\:\:\:$}}{\small{ $\forall\ell\in\mathcal{J}_{k}^{j}\cap\{1,\ldots,m_{1}\}$}} & \tabularnewline
 & \textbf{\small{$\:\:\:\:\:\:\:\:\:$$\:\:\:$}}{\small{ from tracking
neighbors}} & \tabularnewline
{\small{14:}} & \textbf{\small{$\:\:\:\:\:\:$ else }}{\small{Obtain the pdfs}}\textbf{\small{
$\mathbf{\mathcal{F}}_{k,\nu-1}^{\ell},$}} & \tabularnewline
 & \textbf{\small{$\:\:\:\:\:\:\:\:\:$$\:\:\:$}}{\small{ $\forall\ell\in\mathcal{J}_{k}^{j}$
from neighbors}} & \tabularnewline
 & \textbf{\small{$\:\:\:\:\:\:$end if}} & \tabularnewline
{\small{15:}} & \textbf{\small{$\:\:\:\:\:\:$}}{\small{Compute the new pdf $\mathbf{\mathcal{F}}_{k,\nu}^{j}$ }} & \tabularnewline
 & \textbf{\small{$\:\:\:\:\:\:$$\:\:\:$}}{\small{using the LogOP (\ref{log_agreement_equation})}} & \tabularnewline
 & \textbf{\small{$\:\:\:$end for}} & \tabularnewline
{\small{16:}} & \textbf{\small{$\:\:\:$}}{\small{Set $\mathbf{\mathcal{F}}_{k}^{j}=\mathbf{\mathcal{F}}_{k,n_{\textrm{loop}}}^{j}$}} & \tabularnewline
 & \textbf{\small{end while}} & \tabularnewline
\end{tabular}\end{singlespace}
\end{algorithm}

\begin{assumption} The communication network topologies $\mathcal{G}_{k}$
and $\mathcal{D}_{k}$ are SC and the weights $a_{k}^{j\ell}$ are
such that the digraph $\mathcal{D}_{k}$ is balanced. The weights
$a_{k,\nu-1}^{j\ell},\forall j,\ell\in\{1,\ldots,m\}$ and the matrix
$P_{k,\nu-1}$ have the following properties: (i) The weights are
the same for all consensus loops within each time instants, i.e.,
$a_{k,\nu-1}^{j\ell}=a_{k}^{j\ell}$ and $P_{k,\nu-1}=P_{k},\forall\nu\in\mathbb{N}$.
Moreover, $P_{k}$ can be decomposed into four parts $P_{k}=\left[\begin{smallmatrix}P_{k1} & P_{k2}\\
P_{k3} & P_{k4}
\end{smallmatrix}\right]$, where $P_{k1}\in\mathbb{R}^{m_{1}\times m_{1}}$, $P_{k2}=\mathbb{R}^{m_{1}\times m_{2}}$,
$P_{k3}\in\mathbb{R}^{m_{2}\times m_{1}}$, and $P_{k4}\in\mathbb{R}^{m_{2}\times m_{2}}$.
(ii) If $j\in\{1,\ldots,m_{1}\}$, then $a_{k}^{j\ell}>0$ if and
only if $\ell\in\mathcal{J}_{k}^{j}\cap\{1,\ldots,m_{1}\}$, else
$a_{k}^{j\ell}=0$; hence $P_{k2}=\mathbf{0}^{m_{1}\times m_{2}}$.
Moreover, $P_{k1}$ is balanced, i.e., $\sum_{\ell\in\mathcal{J}_{k}^{j}}a_{k}^{j\ell}=\sum_{r|j\in\mathcal{J}_{k}^{r}}a_{k}^{rj}$,
where $j,\ell,r\in\{1,\ldots,m_{1}\}$; (iii) If $j\in\{m_{1}+1,\ldots,m\}$,
then $a_{k}^{j\ell}>0$ if and only if $\ell\in\mathcal{J}_{k}^{j}$,
else $a_{k}^{j\ell}=0$; (iv) The matrix $P_{k}$ is row stochastic,
i.e., $\sum_{\ell=1}^{m}a_{k}^{j\ell}=1$. \label{assump:weights4}
\hfill $\Box$ \end{assumption}

\begin{theorem} \label{thm:Hierarchical-LogOP-SC-balanced-digraphs}
\emph{(Hierarchical Consensus using the LogOP on SC Balanced Digraphs)}
Under Assumptions \ref{assump:nonnegative_pdf} and \ref{assump:weights4},
using the LogOP (\ref{log_agreement_equation}), each $\mathbf{\mathcal{F}}_{k,\nu}^{j}$
globally exponentially converges pointwise to the pdf $\mathbf{\mathcal{F}}_{k}^{\star}$
given by: \vspace{-15pt}

\begin{equation}
\mathbf{\mathcal{F}}_{k}^{\star}=p_{k}^{\star}(\boldsymbol{x}_{k})=\frac{\prod_{i=1}^{m_{1}}\left(p_{k,0}^{i}(\boldsymbol{x}_{k})\right)^{\frac{1}{m_{1}}}}{\int_{\mathcal{X}}\prod_{i=1}^{m_{1}}\left(p_{k,0}^{i}(\boldsymbol{x}_{k})\right)^{\frac{1}{m_{1}}}\: d\mu(\boldsymbol{x}_{k})}\label{WPH_eqn6}
\end{equation}
at a rate faster or equal to $\sqrt{\lambda_{m_{1}-1}(P_{k1}^{T}P_{k1})}=\sigma_{m_{1}-1}(P_{k1})$.
Only the initial estimates of the tracking agents contribute to the
consensual pdf $\mathbf{\mathcal{F}}_{k}^{\star}$. Furthermore, their
induced measures converge in total variation, i.e., $\lim_{\nu\rightarrow\infty}\mu_{\mathbf{\mathcal{F}}_{k,\nu}^{j}}\xrightarrow{\textrm{T.V.}}\mu_{\mathcal{F}_{k}^{\star}},\thinspace\forall j\in\{1,\ldots,m\}$.
\end{theorem}
\begin{IEEEproof}
The matrix $P_{k1}$ conforms to the balanced digraph $\mathcal{D}_{k}$.
Let $\mathbf{1}_{m_{1}}=[1,1,\ldots,1]^{T}$, with $m_{1}$ elements.
Similar to the proof of Theorem \ref{thm:LogOP-SC-balanced-digraphs},
we get $P_{k1}$ is a primitive matrix and $\lim_{\nu\rightarrow\infty}P_{k1}^{\nu}=\frac{1}{m_{1}}\mathbf{1}_{m_{1}}\mathbf{1}_{m_{1}}^{T}$. 

Next, we decompose $\mathbf{\mathcal{U}}_{k,\nu}$ from equation (\ref{log2_PW_eqn2})
into two parts such that $\mathbf{\mathcal{U}}_{k,\nu}=\left[\mathbf{\mathcal{Y}}_{k,\nu};\thinspace\mathcal{Z}_{k,\nu}\right]$,
where $\mathcal{Y}_{k,\nu}=\left(\mathbf{\mathcal{H}}_{k,\nu}^{1},\ldots,\mathbf{\mathcal{H}}_{k,\nu}^{m_{1}}\right)^{T}$
and $\mathbf{\mathcal{Z}}_{k,\nu}=\left(\mathbf{\mathcal{H}}_{k,\nu}^{m_{1}+1},\ldots,\mathbf{\mathcal{H}}_{k,\nu}^{m}\right)^{T}$.
Since $P_{k2}$ is a zero matrix, (\ref{log2_PW_eqn2}) generalizes
and hierarchically decomposes to: \vspace{-20pt}

\begin{eqnarray}
\mathbf{\mathcal{Y}}_{k,\nu+1} & = & P_{k1}^{\nu}\mathbf{\mathcal{Y}}_{k,0},\thinspace\forall\nu\in\mathbb{N}\label{WPH_eqn2_part1}\\
\mathbf{\mathcal{Z}}_{k,\nu+1} & = & P_{k3}\mathbf{\mathcal{Y}}_{k,\nu}+P_{k4}\mathbf{\mathcal{Z}}_{k,\nu},\thinspace\forall\nu\in\mathbb{N}\label{WPH_eqn2_part2}
\end{eqnarray}
Combining equation (\ref{WPH_eqn2_part1}) with the previous result
gives $\lim_{\nu\rightarrow\infty}\mathbf{\mathcal{Y}}_{k,\nu}=\frac{1}{m_{1}}\mathbf{1}_{m_{1}}\mathbf{1}_{m_{1}}^{T}\mathbf{\mathcal{Y}}_{k,0}$.
Thus $\lim_{\nu\rightarrow\infty}\mathbf{\mathcal{H}}_{k,\nu}^{j}=\mathbf{\mathcal{H}}_{k}^{\star}=\frac{1}{m_{1}}\mathbf{1}_{m_{1}}^{T}\mathbf{\mathcal{Y}}_{k,0}=\frac{1}{m_{1}}\sum_{i=1}^{m_{1}}\mathbf{\mathcal{H}}_{k,0}^{i},\forall j\in\{1,\ldots,m_{1}\}$.
From the proof of Theorem \ref{thm:LogOP-SC-balanced-digraphs}, we
get $\mathbf{\mathcal{F}}_{k,\nu}^{j},\forall j\in\{1,\ldots,m_{1}\}$
globally exponentially converges pointwise to $\mathbf{\mathcal{F}}_{k}^{\star}$
given by (\ref{WPH_eqn6}) with a rate faster or equal to $\sigma_{m_{1}-1}(P_{k1})$. 

Since $\mathcal{G}(k)$ is strongly connected, information from the
tracking agents reach the non-tracking agents. Taking the limit of
equation (\ref{WPH_eqn2_part2}) and substituting the above result
gives: \vspace{-15pt}

\begin{equation}
\lim_{\nu\rightarrow\infty}\mathbf{\mathcal{Z}}_{k,\nu+1}=\frac{1}{m_{1}}P_{k3}\mathbf{1}_{m_{1}}\mathbf{1}_{m_{1}}^{T}\mathbf{\mathcal{Y}}_{k,0}+P_{k4}\lim_{\nu\rightarrow\infty}\mathbf{\mathcal{Z}}_{k,\nu}\label{WPH_eqn4}
\end{equation}
Let $\mathbf{1}_{m_{2}}=[1,1,\ldots,1]^{T}$, with $m_{2}$ elements.
Since $P_{k}$ is row stochastic, we get $P_{k3}\mathbf{1}_{m_{1}}=[\mathbf{I}-P_{k4}]\mathbf{1}_{m_{2}}$.
Hence, from equation (\ref{WPH_eqn4}), we get $\lim_{\nu\rightarrow\infty}\mathbf{\mathcal{Z}}_{k,\nu}=\frac{1}{m_{1}}\mathbf{1}_{m_{2}}\mathbf{1}_{m_{1}}^{T}\mathbf{\mathcal{Y}}_{k,0}$.
Moreover, the inessential states die out geometrically fast \cite[pp. 120]{Ref:Seneta06}.
Hence $\lim_{\nu\rightarrow\infty}\mathbf{\mathcal{H}}_{k,\nu}^{j}=\mathbf{\mathcal{H}}_{k}^{\star}=\frac{1}{m_{1}}\mathbf{1}_{m_{1}}^{T}\mathbf{\mathcal{Y}}_{k,0},\forall j\in\{m_{1}+1,\ldots,m\}$.
Hence, the estimates of the non-tracking agents $\mathbf{\mathcal{F}}_{k,\nu}^{j},\forall j\in\{m_{1}+1,\ldots,m\}$
also converge pointwise geometrically fast to the same consensual
pdf $\mathbf{\mathcal{F}}_{k}^{\star}$ given by (\ref{WPH_eqn6}).
By Lemma \ref{lemma:convegence-distribution-measure} we get $\lim_{\nu\rightarrow\infty}\mu_{\mathbf{\mathcal{F}}_{k,\nu}^{j}}\xrightarrow{\textrm{T.V.}}\mu_{\mathbf{\mathcal{F}}_{k}^{\star}},\forall j\in\{1,\ldots,m\}$. 
\end{IEEEproof}
Note that Theorem \ref{thm:BCF-LogOP-nloop} can be directly applied
from Section \ref{sub:Bayesian-Consensus-Filtering} to find the minimum
number of consensus loops $n_{\textrm{loop}}$ for achieving $\epsilon$-convergence
in a given communication network topology or for designing the $P_{k1}$
matrix for a given number of consensus loops. A simulation example
of Hierarchical BCF--LogOP algorithm for tracking orbital debris in
space is discussed in the next section. 

\begin{figure}
\begin{centering}
\includegraphics[bb=0bp 0bp 583bp 454bp,clip,width=3.3in]{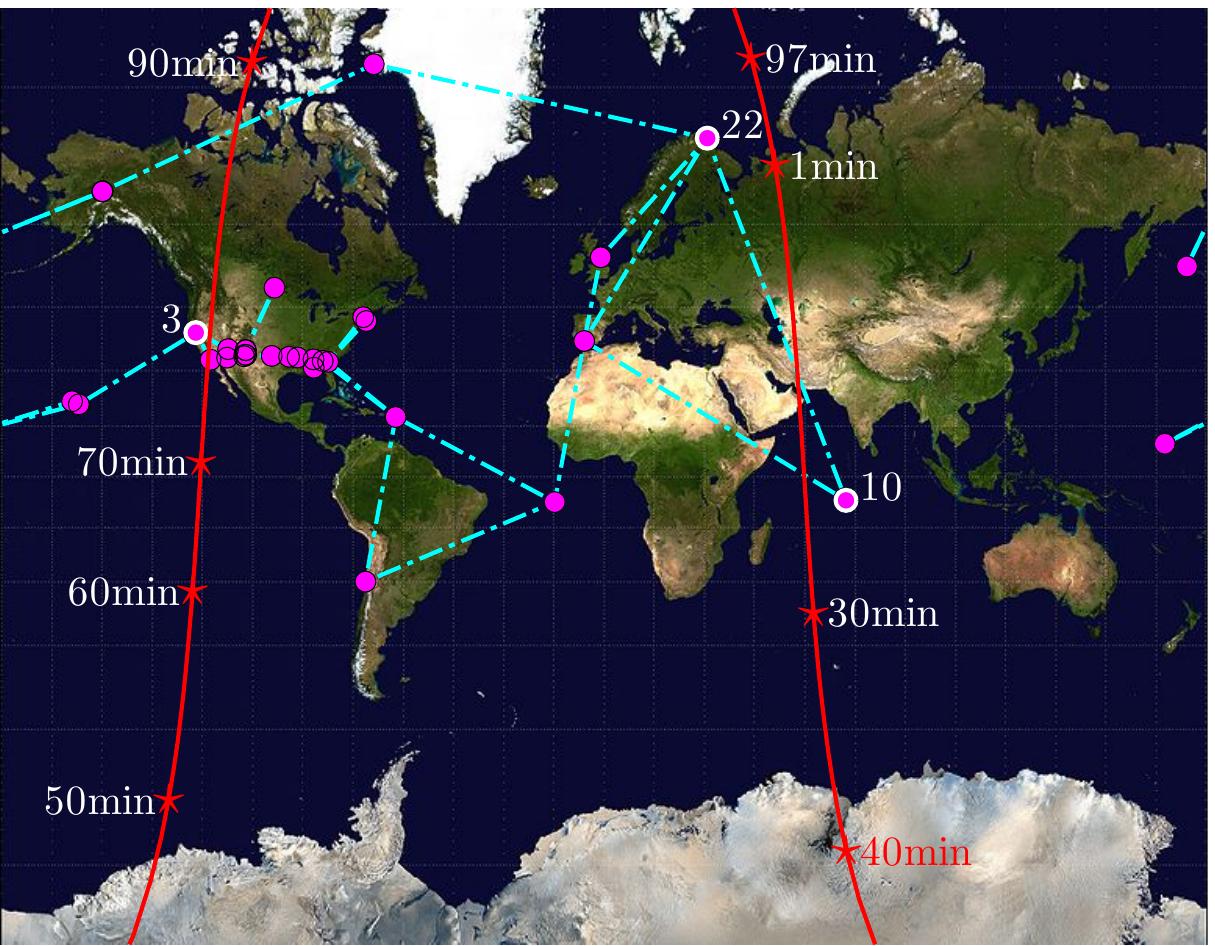}
\par\end{centering}

\centering{}\caption{The SSN locations are shown along with their static SC balanced communication
network topology. The orbit of the Iridium--33 debris is shown in
red, where \textcolor{red}{\large{$\star$}} marks its actual position
during particular time instants. \label{comm_net_topo}}
\end{figure}

\section{Numerical Example \label{sec:Numerical-Example}}

Currently, there are over ten thousand objects in Earth orbit, of
size $0.5$ cm or greater, and almost $95\%$ of them are nonfunctional
space debris. These debris pose a significant threat to functional
spacecraft and satellites in orbit. The US has established the Space
Surveillance Network (SSN) for ground based observations of the orbital
debris using radars and optical telescopes \cite{Ref:Chatters09,Ref:Vallado12}.
In February $2009$, the Iridium--33 satellite collided with the Kosmos--2251
satellite and a large number of debris fragments were created. In
this section, we use the Hierarchical BCF--LogOP Algorithm to track
one of the Iridium--33 debris created in this collision. The orbit
of this debris around Earth and the location of SSN sensors are shown
in Fig. \ref{comm_net_topo}. 

\begin{figure}
\begin{centering}
\begin{tabular}{cc}
\includegraphics[bb=0bp 0bp 333bp 303bp,clip,width=1.45in]{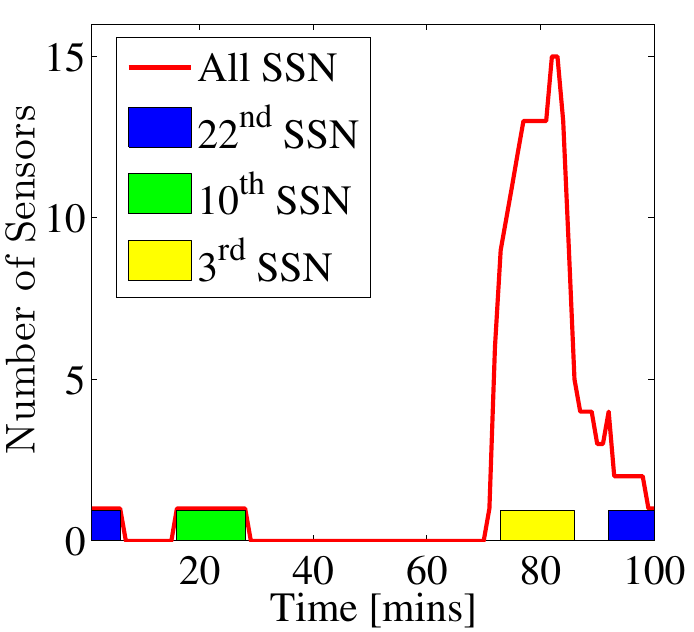} & \includegraphics[bb=0bp 0bp 349bp 303bp,clip,width=1.5in]{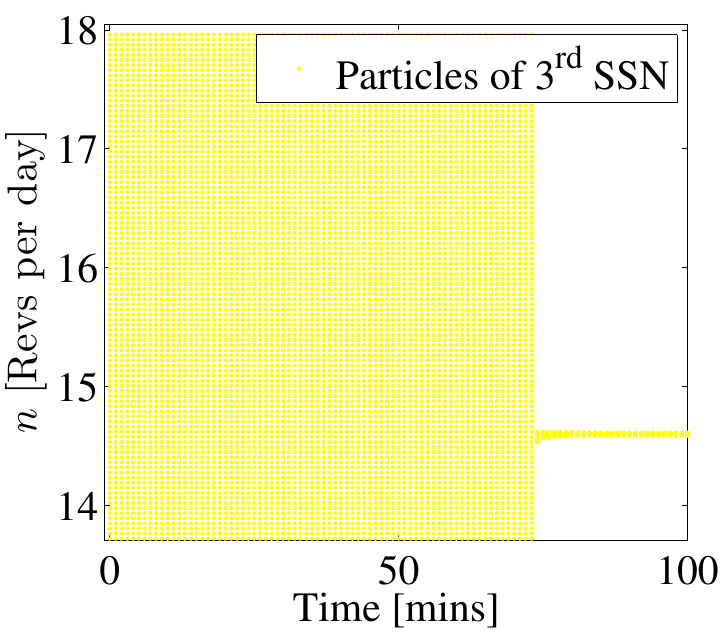}\tabularnewline
(a) & (b)\tabularnewline
\includegraphics[bb=0bp 0bp 349bp 303bp,clip,width=1.5in]{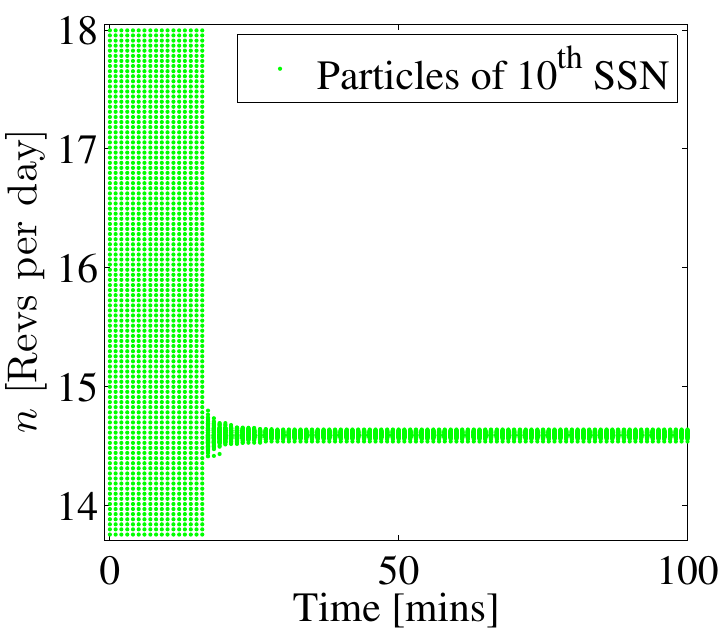} & \includegraphics[bb=0bp 0bp 349bp 303bp,clip,width=1.5in]{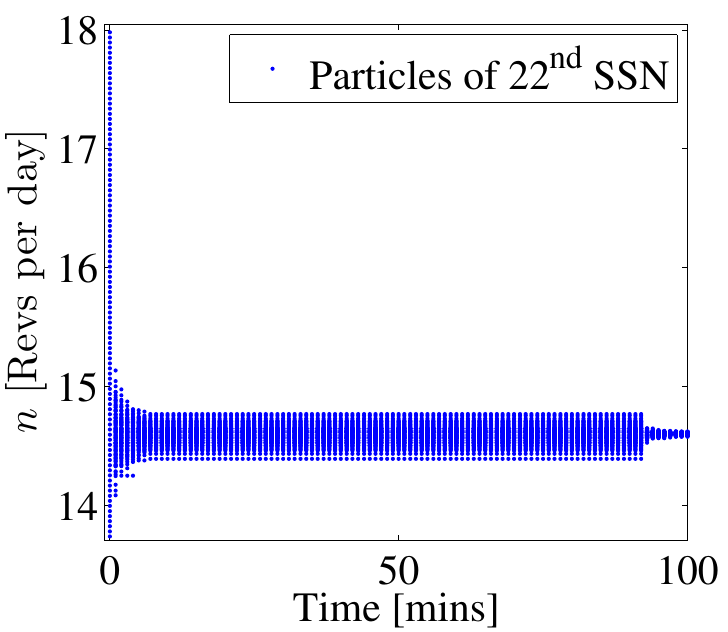}\tabularnewline
(c) & (d)\tabularnewline
\end{tabular}
\par\end{centering}

\caption{(a) Number of SSN sensors observing debris. Trajectories of particles
for stand-alone Bayesian filters for (b) $3^{\textrm{rd}}$, (c) $10^{\textrm{th}}$,
and (d) $22^{\textrm{nd}}$ SSN sensor. \label{fig:SSN_obs_debris_BF}}

\vspace{-10pt}
\end{figure}

The actual two-line element set (TLE) of the Iridium--33 debris was
accessed from North American Aerospace Defense Command (NORAD) on
$4^{\textrm{th}}$ Dec $2013$. The nonlinear Simplified General Perturbations
(SGP4) model, which uses an extensive gravitational model and accounts
for the drag effect on mean motion \cite{Ref:Hoots80,Ref:Vallado08},
is used as the target dynamics model. The communication network topology
of the SSN is assumed to be a static SC balanced graph, as shown in
Figure \ref{comm_net_topo}. If the debris is visible above the sensor's
horizon, then it is assumed to create a single measurement during
each time step of one minute. The heterogeneous measurement model
of the $j^{\textrm{th}}$ sensor is given by: 
\[
\boldsymbol{z}_{k}^{j}=\boldsymbol{x}_{k}+\boldsymbol{w}_{k}^{j},\textrm{ where }\boldsymbol{w}_{k}^{j}=\mathcal{N}\left(0,(1000+50j)\times\mathbf{I}\right),
\]
where $\boldsymbol{x}_{k}\in\mathbb{R}^{3}$ is the actual location
of the debris and the additive Gaussian measurement noise depends
on the sensor number. Since it is not possible to implement the SGP4
target dynamics on distributed estimation algorithms discussed in
the literature \cite{Ref:Borkar82}--\cite{Ref:Hadaegh07}, we compare
the performance of our Hierarchical BCF--LogOP algorithm (\textbf{Algorithm
\ref{alg:HBCF-LogOP}}) against the Hierarchical BCF--LinOP algorithm,
where the LinOP is used during the consensus stage. 

In this simulation example, we simplify the debris tracking problem
by assuming only the mean motion ($n$) of the debris is unknown.
The objective of this simulation example is to estimate $n$ of the
Iridium--33 debris within $100$ minutes. Hence, each sensor knows
the other TLE parameters of the debris and an uniform prior distribution
($\mathbf{\mathcal{F}}_{0}^{j}$) is assumed. Note that at any time
instant, only a few of the SSN sensors can observe the debris, as
shown in Fig \ref{fig:SSN_obs_debris_BF}(a). The results of three
stand-alone Bayesian filters, implemented using particle filters with
resampling \cite{Ref:Arulampalam02}, are shown in Fig \ref{fig:SSN_obs_debris_BF}(b-d).
Note that the estimates of the $22^{\textrm{nd}}$ and $10^{\textrm{th}}$
sensors initially do not converge due to large measurement error,
in spite of observing the debris for some time. The estimates of the
$3^{\textrm{rd}}$ sensor does converge when it is able to observe
the debris after 70 minutes. Hence we propose to use the Hierarchical
BCF--LogOP algorithm where the consensual distribution is updated
as and when sensors observe the debris. 

\begin{figure}
\begin{centering}
\begin{tabular}{cc}
\includegraphics[bb=0bp 0bp 336bp 297bp,clip,width=1.5in]{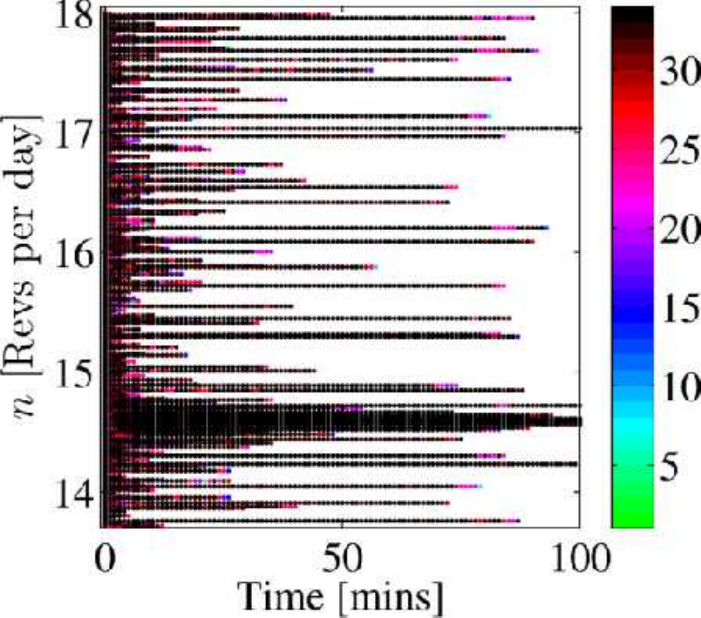} & \includegraphics[bb=0bp 0bp 336bp 297bp,clip,width=1.5in]{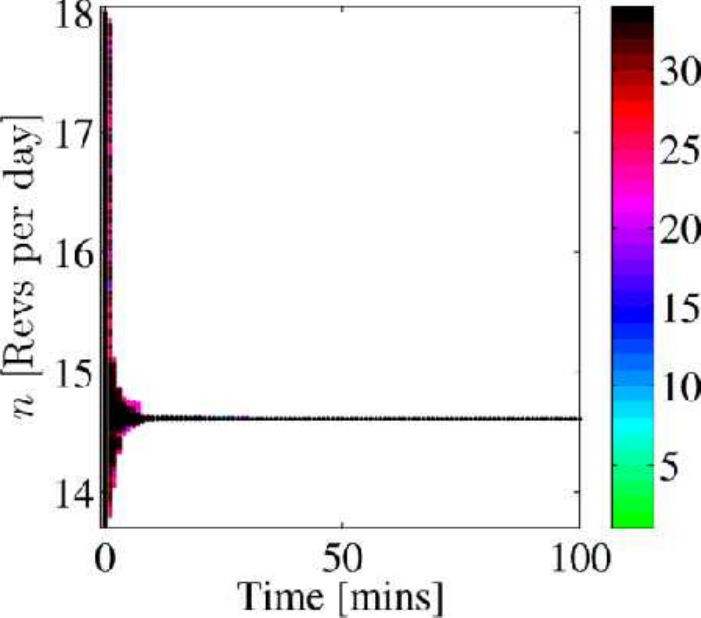}\tabularnewline
(a) & (b)\tabularnewline
\includegraphics[bb=0bp 0bp 400bp 310bp,clip,width=1.5in]{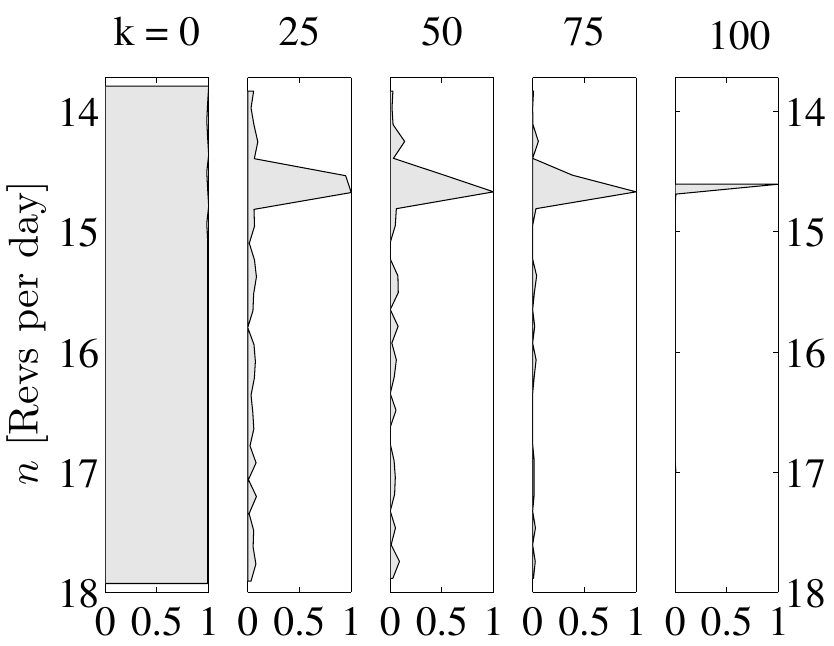} & \includegraphics[bb=0bp 0bp 400bp 310bp,clip,width=1.5in]{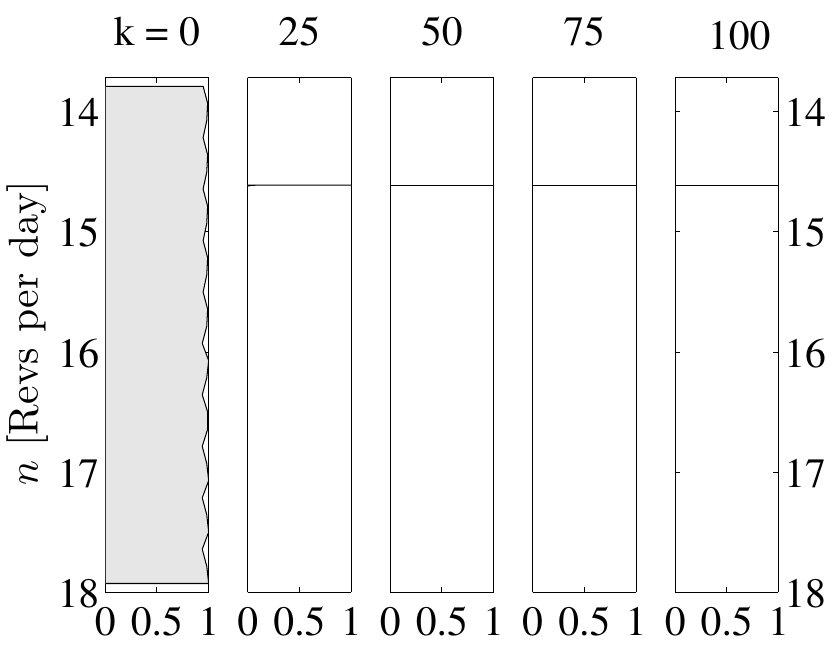}\tabularnewline
(c) & (d)\tabularnewline
\end{tabular}
\par\end{centering}

\caption{Trajectories of particles of all sensors for (a) Hierarchical BCF--LinOP
and (b) Hierarchical BCF--LogOP. The color-bar on the right denotes
the $33$ SSN sensors. Evolution of the consensual probability distribution
for (c) Hierarchical BCF--LinOP and (d) Hierarchical BCF--LogOP. \label{fig:Trajectory_consensus_LinOP_LogOP}}
\end{figure}

Particle filters with resampling are used to evaluate the Bayesian
filters and communicate pdfs in the Hierarchical BCF algorithms. $100$
particles are used by each sensor and $10$ consensus loops are executed
during each time step of one minute. The trajectories of all the particles
of the sensors in the Hierarchical BCF algorithm using LinOP and LogOP
and their respective consensual probability distributions at different
time instants are shown in Figure \ref{fig:SSN_obs_debris_BF}(a-d).
As expected, all the sensors converge on the correct value of $n$
of $14.6$ revs per day. The Hierarchical BCF--LinOP estimates are
multimodal for the first $90$ minutes. On the other hand, the Hierarchical
BCF--LogOP estimates converges to the correct value within the first
10 minutes because the LogOP algorithm efficiently communicates the
best consensual estimate to other sensors during each time step and
achieves consensus across the network.

\section{Conclusion \label{sec:Conclusion}}

In this paper, we extended the scope of distributed estimation algorithms
in a Bayesian filtering framework in order to simultaneously track
targets, with general nonlinear time-varying target dynamic models,
using a strongly connected network of heterogeneous agents, with general
nonlinear time-varying measurement models. We introduced the Bayesian
filter with/without measurement exchange to generate local estimated
pdfs of the target's states. We compared the LinOP and LogOP methods
of combining local posterior pdfs and determined that LogOP is the
superior scheme. The LogOP algorithm on SC balanced digraph converges
globally exponentially, and the consensual pdf minimizes the information
lost during the consensus stage because it minimizes the sum of KL
divergences to each locally estimated probability distribution. We
also explored several methods of communicating pdfs across the sensor
network. We introduced the BCF algorithm, where the local estimated
posterior pdfs of the target's states are first updated using the
Bayesian filter and then recursively combined during the consensus
stage using LogOP, so that the agents can track a moving target and
also maintain consensus across the network. Conditions for exponential
convergence of the BCF algorithm and constraints on the communication
network topology have been studied. The Hierarchical BCF algorithm,
where some of the agents do not observe the target, has also been
investigated. Simulation results demonstrate the effectiveness of
the BCF algorithms for nonlinear distributed estimation problems.

\section*{Acknowledgment}

The authors would like to thank F. Hadaegh, D. Bayard, S. Hutchinson,
P. Voulgaris, M. Egerstedt, A. Gupta, A. Dani, D. Morgan, S. Sengupta,
and A. Olshevsky for stimulating discussions about this paper.

\appendices{}

\section{Proof of Theorem \ref{thm:LinOP-SC-digraphs} \label{sec:Proof-of-Theorem-BCF-LinOP-Static}}

Under Assumption \ref{assump:weights1}, $P_{k}$ is a nonnegative,
row stochastic and irreducible matrix. Similar to the proof in \cite{Ref:Saber07},
all the diagonal entries of $P_{k}$ are positive, then $P_{k}^{m-1}>0$
and \cite[Theorem 8.5.2, pp. 516]{Ref:Horn85} implies that $P_{k}$
is a primitive matrix. Since $P_{k}$ is a regular matrix, it has
only one recurrent class which is aperiodic \cite[pp. 127]{Ref:Seneta06}.
Since $P_{k}$ is row stochastic, $\mathbf{1}$ is its right eigenvector
corresponding to the eigenvalue $1$, i.e., $P_{k}\mathbf{1}=1\mathbf{1}$.
Moreover, according to the Gershgorin Disc Theorem (cf. \cite[pp. 344]{Ref:Horn85}),
all the eigenvalues of $P_{k}$ are located in the unit circle, i.e.
the spectral radius $\rho(P_{k})=1$. Hence Perron--Frobenius theorem
(cf. \cite[pp. 3]{Ref:Seneta06}) states that there exists a left
eigenvector $\boldsymbol{\pi}$ of $P_{k}$ corresponding to the eigenvalue
$1$ which is unique to constant multiples, i.e., $P_{k}^{T}\boldsymbol{\pi}=1\boldsymbol{\pi}$.
The ergodic theorem for primitive Markov chains (cf. \cite[pp. 119]{Ref:Seneta06})
states that $P_{k}$ has an unique stationary distribution given by
the solution of the normalizing condition $\boldsymbol{\pi}^{T}\mathbf{1}=1$
and $\lim_{\nu\rightarrow\infty}P_{k}^{\nu}=\mathbf{1}\boldsymbol{\pi}^{T}$.

Obviously, (\ref{PW_eqn}) generalizes to $\mathbf{\mathcal{W}}_{k,\nu}=P_{k}^{\nu}\mathbf{\mathcal{W}}_{k,0},\forall\nu\in\mathbb{N}$,
where $\mathbf{\mathcal{W}}_{k,0}=\left(\mathbf{\mathcal{F}}_{k,0}^{1},\ldots,\mathbf{\mathcal{F}}_{k,0}^{m}\right)^{T}$.
Hence we get $\lim_{\nu\rightarrow\infty}\mathbf{\mathcal{W}}_{k,\nu}=\mathbf{1}\boldsymbol{\pi}^{T}\mathbf{\mathcal{W}}_{k,0}$.
Thus, each $\mathbf{\mathcal{F}}_{k,\nu}^{j}$ converges pointwise
to the consensual pdf $\mathbf{\mathcal{F}}_{k}^{\star}=\boldsymbol{\pi}^{T}\mathbf{\mathcal{W}}_{k,0}$.
Note that $\sum_{i=1}^{m}\pi_{i}=1$, as needed for $\mathbf{\mathcal{F}}_{k}^{\star}$
to be a valid probability distribution, where $\pi_{i}$ is the individual
element of the vector $\boldsymbol{\pi}$. By Lemma \ref{lemma:convegence-distribution-measure},
the measure induced by $\mathbf{\mathcal{F}}_{k,\nu}^{j}$ on $\mathscr{X}$
converges in total variation to the measure induced by $\mathbf{\mathcal{F}}_{k}^{\star}$
on $\mathscr{X}$, i.e., $\lim_{\nu\rightarrow\infty}\mu_{\mathbf{\mathcal{F}}_{k,\nu}^{j}}\xrightarrow{\textrm{T.V.}}\mu_{\mathbf{\mathcal{F}}_{k}^{\star}}$.
\hfill $\blacksquare$ 

\bibliographystyle{ieeetr}
\bibliography{IEEEabrv,SapBib}

\begin{thebibliography}{10}

\bibitem{Ref:Saber04}
R.~Olfati-Saber and R.~Murray, ``Consensus problems in networks of agents with
  switching topology and time-delays,'' {\em {IEEE} Trans. Autom. Control},
  vol.~49, no.~9, pp.~1520 -- 1533, 2004.

\bibitem{Ref:Jadbabaie03}
A.~Jadbabaie, J.~Lin, and A.~S. Morse, ``Coordination of groups of mobile
  autonomous agents using nearest neighbor rules,'' {\em {IEEE} Trans. Autom.
  Control}, vol.~48, no.~6, pp.~988 -- 1001, 2003.

\bibitem{Ref:Shah06}
S.~Boyd, A.~Ghosh, B.~Prabhakar, and D.~Shah, ``Randomized gossip algorithms,''
  {\em {IEEE} Trans. Inf. Theory}, vol.~52, pp.~2508 -- 2530, June 2006.

\bibitem{Ref:Boyd04}
L.~Xiao and S.~Boyd, ``Fast linear iterations for distributed averaging,'' {\em
  Syst. Control Lett.}, vol.~53, pp.~65 -- 78, 2004.

\bibitem{Ref:Ren05TAC}
W.~Ren and R.~W. Beard, ``Consensus seeking in multiagent systems under
  dynamically changing interaction topologies,'' {\em {IEEE} Trans. Autom.
  Control}, vol.~50, pp.~655 -- 661, May 2005.

\bibitem{Ref:Saber07}
R.~Olfati-Saber, J.~A. Fax, and R.~M. Murray, ``Consensus and cooperation in
  networked multi-agent systems,'' {\em Proc. of the IEEE}, vol.~95, no.~1,
  pp.~215--233, 2007.

\bibitem{Ref:Tsitsiklis05}
V.~D. Blondel, J.~M. Hendrickx, A.~Olshevsky, and J.~N. Tsitsiklis,
  ``Convergence in multiagent coordination, consensus, and flocking,'' in {\em
  44th IEEE Conf. Decision Control}, (Seville, Spain), Dec. 2005.

\bibitem{Ref:Tsitsiklis86}
J.~Tsitsiklis, D.~Bertsekas, and M.~Athans, ``Distributed asynchronous
  deterministic and stochastic gradient optimization algorithms,'' {\em {IEEE}
  Trans. Autom. Control}, vol.~31, no.~9, pp.~803 -- 812, 1986.

\bibitem{Ref:Martinez10}
M.~Zhu and S.~Mart\'{i}nez, ``On distributed optimization under inequality and
  equality constraints via penalty primal-dual methods,'' in {\em Amer. Control
  Conf.}, (Baltimore, US), pp.~2434--2439, 2010.

\bibitem{Ref:Nedic09}
A.~Nedic and A.~Ozdaglar, {\em Convex Optimization in Signal Processing and
  Communications}, ch.~Cooperative distributed multi-agent optimization,
  pp.~340 -- 386.
\newblock Cambridge University Press, 2009.

\bibitem{Ref:Borkar82}
V.~Borkar and P.~Varaiya, ``Asymptotic agreement in distributed estimation,''
  {\em {IEEE} Trans. Autom. Control}, vol.~27, no.~3, pp.~650 -- 655, 1982.

\bibitem{Ref:Moura11}
S.~Kar and J.~M.~F. Moura, ``Convergence rate analysis of distributed gossip
  (linear parameter) estimation: Fundamental limits and tradeoffs,'' {\em
  {IEEE} J. Sel. Topics Signal Process.}, vol.~5, pp.~674 -- 690, Aug 2011.

\bibitem{Ref:Giannakis09}
I.~Schizas, G.~Mateos, and G.~Giannakis, ``Distributed {LMS} for
  consensus-based in-network adaptive processing,'' {\em {IEEE} Trans. Signal
  Process.}, vol.~57, pp.~2365 -- 2382, June 2009.

\bibitem{Ref:Murray05}
D.~Spanos, R.~Olfati-Saber, and R.~M. Murray, ``Distributed sensor fusion using
  dynamic consensus,'' in {\em Proc. {IFAC}}, 2005.

\bibitem{Ref:Coates04}
M.~Coates, ``Distributed particle filters for sensor networks,'' in {\em Proc.
  3rd Int. Symp. Inform. Process. Sensor Networks}, (New York, USA),
  pp.~99--107, 2004.

\bibitem{Ref:Leonard10}
F.~Zhang and N.~E. Leonard, ``Cooperative filters and control for cooperative
  exploration,'' {\em {IEEE} Trans. Autom. Control}, vol.~55, no.~3,
  pp.~650--663, 2010.

\bibitem{Ref:Ahmed13}
N.~Ahmed, J.~Schoenberg, and M.~Campbell, {\em Robotics: Science and Systems
  VIII}, ch.~Fast Weighted Exponential Product Rules for Robust General
  Multi-Robot Data Fusion, pp.~9--16.
\newblock MIT Press, 2013.

\bibitem{Ref:Demetriou11}
M.~Demetriou and D.~Uci\'{n}ski, ``State estimation of spatially distributed
  processes using mobile sensing agents,'' in {\em Amer. Control Conf.}, (San
  Francisco, CA, USA), pp.~1770--1776, 2011.

\bibitem{Ref:Moura10}
A.~Dimakis, S.~Kar, J.~Moura, M.~Rabbat, and A.~Scaglione, ``Gossip algorithms
  for distributed signal processing,'' {\em Proc. of the {IEEE}}, vol.~98,
  no.~11, pp.~1847--1864, 2010.

\bibitem{Ref:Hadaegh07}
B.~A\c{c}ikme\c{s}e, F.~Y. Hadaegh, D.~P. Scharf, and S.~R. Ploen,
  ``Formulation and analysis of stability for spacecraft formations,'' {\em
  {IET} Control Theory Appl.}, vol.~1, no.~2, pp.~461--474, 2007.

\bibitem{Ref:Gupta13}
Y.~Xu, V.~Gupta, and C.~Fischione, ``Distributed estimation,'' in {\em
  E-reference Signal Processing}, Elsevier, 2013.
\newblock Editors: Rama Chellappa, Sergios Theodoridis.

\bibitem{Ref:Saber09}
R.~Olfati-Saber, ``Kalman-consensus filter : Optimality, stability, and
  performance,'' in {\em 48th IEEE Conf. Decision Control}, (Shanghai, China),
  pp.~7036--7042, December 2009.

\bibitem{Ref:DeGroot75}
M.~H. DeGroot, {\em Probability and Statistics}.
\newblock Cambridge, Massachusetts: Addison-Wesley, 1975.

\bibitem{Ref:Jeffreys61}
H.~Jeffreys, {\em Theory of Probability}.
\newblock Oxford: Clarendon Press, 1961.

\bibitem{Ref:Subrahmaniam79}
K.~Subrahmaniam, {\em A Primer in Probability}.
\newblock New York, NY: M. Dekker, 1979.

\bibitem{Ref:Pearl88}
J.~Pearl, {\em Probabilistic Reasoning in Intelligent Systems: Networks of
  Plausible Inference}.
\newblock San Mateo, CA: Morgan Kaufmann, 1988.

\bibitem{Ref:Kalman60}
R.~E. Kalman, ``A new approach to linear filtering and prediction problems,''
  {\em Trans. {ASME} J. Basic Eng.}, vol.~82, no.~Series D, pp.~35--45, 1960.

\bibitem{Ref:Julier04}
S.~J. Julier and J.~K. Uhlmann, ``Unscented filtering and nonlinear
  estimation,'' {\em Proc. of the {IEEE}}, vol.~92, no.~3, pp.~401--422, 2004.

\bibitem{Ref:Thrun00}
S.~Thrun, ``Probabilistic algorithms in robotics,'' {\em AI Magazine}, vol.~21,
  no.~4, pp.~93--109, 2000.

\bibitem{Ref:Burgard96}
W.~Burgard, D.~Fox, D.~Hennig, and T.~Schmidt, ``Estimating the absolute
  position of a mobile robot using position probability grids,'' in {\em Proc.
  of the 14th Nat. Conf. Artificial Intell.}, August 1996.

\bibitem{Ref:Diard03}
J.~Diard, P.~Bessi\`{e}re, and E.~Mazer, ``A survey of probabilistic models,
  using the {B}ayesian programming methodology as a unifying framework,'' in
  {\em 2nd Int. Conf. Computational Intell., Robotics and Autonomous Syst.},
  (Singapore), Dec. 2003.

\bibitem{Ref:Thrun05}
S.~Thrun, W.~Burgard, and D.~Fox, {\em Probabilistic Robotics}.
\newblock Cambridge, Massachusetts: The MIT Press, 2005.

\bibitem{Ref:Boutilier99}
C.~Boutilier, T.~Dean, and S.~Hanks, ``Decision-theoretic planning: Structural
  assumptions and computational leverage,'' {\em J. Artificial Intell.
  Research}, vol.~11, pp.~1--94, 1999.

\bibitem{Ref:Kaelbling98}
L.~P. Kaelbling, M.~L. Littman, and A.~R. Cassandra, ``Planning and acting in
  partially observable stochastic domains,'' {\em Artificial Intell.},
  vol.~101, pp.~99--134, May 1998.

\bibitem{Ref:Punska99}
O.~Punska, ``Bayesian approaches to multi-sensor data fusion,'' Master's
  thesis, Dept. of Eng., Univ. of Cambridge, 1999.

\bibitem{Ref:Arulampalam02}
M.~S. Arulampalam, S.~Maskell, N.~Gordon, and T.~Clapp, ``A tutorial on
  particle filters for online nonlinear/non-{G}aussian {B}ayesian tracking,''
  {\em {IEEE} Trans. Signal Process.}, vol.~50, pp.~174--188, February 2002.

\bibitem{Ref:Lebeltel04}
O.~Lebeltel, P.~Bessiere, J.~Diard, and E.~Mazer, ``Bayesian robot
  programming,'' {\em Autonomous Robots}, vol.~16, pp.~49--79, January 2004.

\bibitem{Ref:Chen05}
M.-H. Chen, ``Bayesian computation: From posterior densities to {B}ayes
  factors, marginal likelihoods, and posterior model probabilities,'' in {\em
  Bayesian Thinking, Modeling and Computation} (D.~K. Dey and C.~R. Rao, eds.),
  Handbook of Statistics, ch.~15, pp.~437 -- 457, Amsterdam: Elsevier, 2005.

\bibitem{Ref:DeGroot74}
M.~H. DeGroot, ``Reaching a consensus,'' {\em J. Amer. Statistical Assoc.},
  vol.~69, no.~345, pp.~688 -- 704, 1960.

\bibitem{Ref:Genest86}
C.~Genest and J.~V. Zidek, ``Combining probability distributions: A critique
  and an annotated bibliography,'' {\em Statistical Sci.}, vol.~1, no.~1,
  pp.~114 -- 135, 1986.

\bibitem{Ref:Bacharach79}
M.~Bacharach, ``Normal {B}ayesian dialogues,'' {\em J. Amer. Statistical
  Assoc.}, vol.~74, no.~368, pp.~837 -- 846, 1979.

\bibitem{Ref:Chatterjee77}
S.~Chatterjee and E.~Seneta, ``Towards consensus: Some convergence theorems on
  repeated averaging,'' {\em J. Appl. Probability}, vol.~14, no.~1, pp.~89 --
  97, 1977.

\bibitem{Ref:French81}
S.~French, ``Consensus of opinion,'' {\em European J. Operational Research},
  vol.~7, pp.~332 -- 340, 1981.

\bibitem{Ref:Velagapudi07}
P.~Velagapudi, O.~Prokopyev, K.~Sycara, and P.~Scerri, ``Maintaining shared
  belief in a large multiagent team,'' in {\em Proc. of {FUSION}}, 2007.

\bibitem{Ref:Yuksel09}
S.~Y\"{u}ksel, ``Stochastic nestedness and the belief sharing information
  pattern,'' {\em {IEEE} Trans. Autom. Control}, vol.~54, no.~12,
  pp.~2773--2786, 2009.

\bibitem{Ref:Fraser12}
C.~S.~R. Fraser, L.~F. Bertuccelli, H.-L. Choi, and J.~P. How, ``A
  hyperparameter consensus method for agreement under uncertainty,'' {\em
  Automatica}, vol.~48, no.~2, pp.~374 -- 380, 2012.

\bibitem{Ref:Hlinka12}
O.~Hlinka, O.~Slu\u{c}iak, F.~Hlawatsch, P.~M. Djuric, and M.~Rupp,
  ``Likelihood consensus and its application to distributed particle
  filtering,'' {\em {IEEE} Trans. Signal Process.}, vol.~60, no.~8,
  pp.~4334--4349, 2012.

\bibitem{Ref:Pavlin10}
G.~Pavlin, P.~Oude, M.~Maris, J.~Nunnink, and T.~Hood, ``A multi-agent systems
  approach to distributed {B}ayesian information fusion,'' {\em Inform.
  Fusion}, vol.~11, pp.~267--282, 2010.

\bibitem{Ref:Jadbabaie12}
A.~Jadbabaie, P.~Molavi, A.~Sandroni, and A.~Tahbaz-Salehi, ``Non-{B}ayesian
  social learning,'' {\em Games and Economic Behavior}, vol.~76, pp.~210--225,
  2012.

\bibitem{Ref:Ross80}
K.~A. Ross, {\em Elementary Analysis: The Theory of Calculus}.
\newblock Springer, 1980.

\bibitem{Ref:Daum05}
F.~Daum, ``Nonlinear filters: beyond the {K}alman filter,'' in {\em {IEEE}
  Aerospace Electron. Syst. Mag.}, vol.~20, pp.~57--69, 2005.

\bibitem{Ref:Cover91}
T.~M. Cover and J.~A. Thomas, {\em Elements of Information Theory}.
\newblock New York, NY: Wiley, 1991.

\bibitem{Ref:Nehorai07}
T.~Zhao and A.~Nehorai, ``Distributed sequential {B}ayesian estimation of a
  diffusive source in wireless sensor networks,'' {\em {IEEE} Trans. Signal
  Process.}, vol.~55, no.~4, pp.~1511 -- 1524, 2007.

\bibitem{Ref:Durrett05}
R.~Durrett, {\em Probability: Theory and Examples}.
\newblock Thomson Brooks, 2005.

\bibitem{Ref:Rus12}
B.~J. Julian, M.~Angermann, M.~Schwager, and D.~Rus, ``Distributed robotic
  sensor networks: An information-theoretic approach,'' {\em Inter. J. of
  Robotics Research}, vol.~31, no.~10, pp.~1134--1154, 2012.

\bibitem{Ref:Weerahandi83}
S.~Weerahandi and J.~V. Zidek, ``Elements of multi-{B}ayesian decision
  theory,'' {\em Ann. of Stat.}, vol.~11, no.~4, pp.~1032 -- 1046, 1983.

\bibitem{Ref:Nash50}
J.~F. Nash, Jr, ``The bargaining problem,'' {\em Econometrica}, vol.~18, no.~2,
  pp.~155 -- 162, 1950.

\bibitem{Ref:Rufo12}
M.~J. Rufo, J.~Mart\'{i}n, and C.~J. P\'{e}rez, ``Log-linear pool to combine
  prior distributions: A suggestion for a calibration-based approach,'' {\em
  Bayesian Analysis}, vol.~7, no.~2, pp.~411--438, 2012.

\bibitem{Ref:Smith05}
A.~Smith, T.~Cohn, and M.~Osborne, ``Logarithmic opinion pools for conditional
  random fields,'' in {\em Proc. Assoc. Computational Linguistics}, (Ann Arbor,
  Michigan), pp.~18--25, 2005.

\bibitem{Ref:Gilardoni93}
G.~L. Gilardoni and M.~K. Clayton, ``On reaching a consensus using
  {D}e{G}root's iterative pooling,'' {\em Ann. Stat.}, vol.~21, no.~1, pp.~391
  -- 401, 1993.

\bibitem{Ref:Chung12}
S.-J. Chung, S.~Bandyopadhyay, I.~Chang, and F.~Y. Hadaegh, ``Phase
  synchronization control of complex networks of {L}agrangian systems on
  adaptive digraphs,'' {\em Automatica}, vol.~49, pp.~1148--1161, May 2013.

\bibitem{Ref:Anderson05}
B.~D.~O. Anderson and J.~B. Moore, {\em Optimal Filtering}.
\newblock Mineola, New York: Dover Publications, 2005.

\bibitem{Ref:Reynolds08}
D.~A. Reynolds, ``Gaussian mixture models,'' {\em Encyclopedia of Biometric
  Recognition}, February 2008.

\bibitem{Ref:McLachlan88}
G.~J. McLachlan and K.~E. Basford, {\em Mixture models : inference and
  applications to clustering}.
\newblock New York, N.Y.: M. Dekker, 1988.

\bibitem{Ref:Kotecha03}
J.~H. Kotecha and P.~M. Djuric, ``Gaussian sum particle filtering,'' {\em
  {IEEE} Trans. Signal Process.}, vol.~51, pp.~2602--2612, Oct. 2003.

\bibitem{Ref:Kramer07}
G.~Kramer and S.~A. Savari, ``Communicating probability distributions,'' {\em
  {IEEE} Trans. Inf. Theory}, vol.~53, pp.~518--525, February 2007.

\bibitem{Ref:Sayed12}
S.-Y. Tu and A.~H. Sayed, ``Diffusion strategies outperform consensus
  strategies for distributed estimation over adaptive networks,'' {\em {IEEE}
  Trans. Signal Process.}, vol.~60, pp.~6217--6234, Dec. 2012.

\bibitem{Ref:Seneta06}
E.~Seneta, {\em Non-negative Matrices and Markov Chains}.
\newblock New York, NY: Springer-Verlag, 2006.

\bibitem{Ref:Chatters09}
M.~E.~P. Chatters and M.~B.~J. Crothers, {\em AU-18 Space Primer}, ch.~Space
  Surveillance Network, pp.~249--258.
\newblock Air University Press, Maxwell Air Force Base, Alabama, 2009.

\bibitem{Ref:Vallado12}
D.~A. Vallado and J.~D. Griesbach, ``Simulating space surveillance networks,''
  in {\em AAS/AIAA Astrodynamics Specialist Conf.}, (Girdwood), 2012.
\newblock Paper AAS 11-580.

\bibitem{Ref:Hoots80}
F.~R. Hoots and R.~L. Roehrich, ``Spacetrack report number 3: Models for
  propagation of {NORAD} element sets,'' tech. rep., U.S. Air Force Aerospace
  Defense Command, Colorado Springs, CO., 1980.

\bibitem{Ref:Vallado08}
D.~Vallado and P.~Crawford, ``{SGP4} orbit determination,'' in {\em AIAA/AAS
  Astrodynamics Specialist Conf.}, 2008.

\bibitem{Ref:Horn85}
R.~A. Horn and C.~R. Johnson, {\em Matrix Analysis}.
\newblock Cambridge, England: Cambridge University Press, 1985.

\end{thebibliography}

\end{document}